\documentclass[11pt]{article}

\usepackage[a4paper,body={160mm,235mm},centering]{geometry}

\usepackage{comment}

\usepackage{amsfonts,amsmath,amsthm,amssymb}

\newtheorem{theorem}{Theorem}[section]
\newtheorem{cor}[theorem]{Corollary}
\newtheorem{lemma}[theorem]{Lemma}
\newtheorem{proposition}[theorem]{Proposition}

\theoremstyle{definition}
\newtheorem{definition}[theorem]{Definition}
\newtheorem{example}[theorem]{Example}
\newtheorem{hypothesis}{Hypothesis}

\theoremstyle{remark}
\newtheorem{remark}[theorem]{Remark}

\makeatletter \@addtoreset{equation}{section}

\makeatother

\def\ds{\begin{displaystyle}}
\def\eds{\end{displaystyle}}
\def\dis{\displaystyle }
\def\<{\langle }
\def\>{\rangle }

\def\R{\mathbb R}
\def\N{\mathbb N}

\def\E{\mathbb E}

\def\P{\mathbb P}

\def\t{\tau}
\def\l{\lambda}
\def\s{\sigma}

\def\calb{{\cal B}}

\def\calf{{\cal F}}
\def\calg{{\cal G}}

\def\calk{{\cal K}}
\def\caln{{\cal N}}
\def\calp{{\cal P}}

\def\call{{\cal L}}
\def\ty{{\widetilde Y}}

\def\tz{{\widetilde Z}}

\newcommand{\rset}{\mathbb{R}}

\newcommand{\nl}{\nolimits}
\newcommand{\ep}{\varepsilon}

\newcommand{\fl}{\longrightarrow}
\newcommand{\lp}{L}
\newcommand{\m}{\mathcal}
\newcommand{\dsp}[1][p]{\mathcal{X}^{#1}}

\title{\bf Differentiability of backward stochastic differential equations in Hilbert spaces with monotone generators.}

\author{Philippe Briand\\[.3em]
\normalsize\it
IRMAR, Universit\'{e} Rennes 1, 35042 Rennes Cedex, FRANCE\\
\normalsize\tt
philippe.briand@univ-rennes1.fr
\and 
Fulvia Confortola\\[.3em]
\normalsize\it
Dipartimento di Matematica, Politecnico di Milano\\
\normalsize\it
piazza Leonardo da Vinci 32, 20133 Milano, Italy\\
\normalsize\tt
confortola@mate.polimi.it}

\date{March 17, 2006}

\begin{document}

\maketitle

\begin{abstract}
The aim of the present paper is to study the regularity properties of the solution of a backward stochastic differential equation with a monotone generator in infinite dimension. We show some applications to the nonlinear Kolmogorov equation and to stochastic optimal control. 
\end{abstract}

\section{Introduction}

In this paper we are concerned mainly with  backward stochastic differential equations (BSDEs for short in the remaining of the paper) in infinite dimension in the markovian framework.  More precisely, we consider the following backward stochastic evolution equation:
\begin{equation}
\label{backprimo}
 \left\{\begin{array}{l}\dis
 dY_\tau=  Z_\tau\, dW_\tau -
  \psi (\tau ,X_\tau, Y_\tau,Z_\tau) \,d\tau,
\quad \tau\in [t,T],
\\\dis
Y_T= \phi(X_T),
\end{array}
\right.
\end{equation}
where $W$ is a cylindrical Wiener process in some Hilbert space $\Xi$. The unknowns are the processes $Y$ and $Z$; $Y$ takes its values in  a Hilbert space $K$ and $Z$ belongs to $L_2(\Xi,K)$, the space of  Hilbert-Schmidt operators from  $\Xi$ to $K$.
The process $X$ takes its values in a Hilbert space $H$ and is the solution of the forward equation
\begin{equation}
\label{stoc}
  \left\{\begin{array}{l}\dis
dX_\tau=AX_\tau \, d\tau+
F(\tau,X_\tau)\, d\tau
+G(\tau,X_\tau)\, dW_\tau,\quad \tau\in [t,T],
\\
\dis X_t=x\in H,
\end{array}\right.
\end{equation}
where $A$ is the generator of a strongly continuous semigroup of
bounded linear operators $\{e^{tA}\}$ in $H$, $F$ and $G$ are
functions with values in $H$ and $L(\Xi,H)$ respectively,
satisfying appropriate Lipschitz conditions.

The above equations (\ref{backprimo}) and (\ref{stoc}) form a so called forward-backward system. 
Under suitable
assumptions on $F$, $G$ and  $\psi$, there exists a unique adapted process
$(X,Y,Z)$ in the space $H\times K\times L_2(\Xi,K)$ solution to this system. The
processes $X,Y,Z$ depend on the values of $x$ and $t$, occurring
as initial conditions in (\ref{stoc}): we may denote them by
$X(\tau, t,x)$, $Y(\tau,t,x)$, $Z(\tau,t,x)$, $\tau\in [t,T]$.

The goal of this work is the study of the continuity and differentiability with respect to the parameters $t$ and $x$ of the process $Y$ solution to the BSDE (\ref{backprimo}) when the coefficient $\psi$ in (\ref{backprimo}) is a monotone operator (see Theorem \ref{ydifferenziabilenelsistema}). 

We show, moreover, some interesting applications of this result to the solvability of nonlinear stochastic Kolmogorov equations 
and to stochastic optimal control problems.

BSDEs in finite and infinite dimensions have been
intensively studied in recent years, starting from the paper by E.
Pardoux and S. Peng \cite{PaPe}: we refer the reader to \cite{ElPQ}, \cite{ElK} and \cite{P3}
  for an exposition
of this subject and to \cite{MaYo} for coupled forward-backward
systems.

The problem of regular dependence of the solution of a stochastic forward-backward system has been studied
in finite dimension by Pardoux, Peng \cite{PaPe2} and by El Karoui, Peng and Quenez \cite{ElPQ},
and, in infinite dimension, by Fuhrman and Tessitore  in \cite{FT1}, \cite{FT2}. In both cases, $\psi$ is assumed to be Lipschitz continuous with respect to $y$ and $z$.
We will assume on $\psi$ a weaker condition: $\psi$ will be Lipschitz continuous only with respect to $z$ and monotone with respect to $y$ in the spirit of the works \cite{Pe1}, \cite{P2} and more recently \cite{BDHPS}.

We should point out that, since we are working in infinite dimensional spaces,  the derivatives are understood in the G\^ateaux sense: for instance
Nemytskii (evaluation) operators are not Fr\'echet differentiable, except in trivial cases.

\smallskip

This result, beside its intrinsic interest in the framework of the theory of BSDEs, allows us to treat nonlinear partial differential equations as the well known Kolmogorov equation: for $t\in[0,T]$ and $x\in H$,
\begin{equation}
\label{kolmogorovnonlineare-intro}
\partial_t v(t,x)+\call_t [v(t,\cdot)](x) +
\psi (t, x,v(t,x),G(t,x)^*\nabla_xv(t,x)) = 0,\qquad  v(T,x)=\phi(x),
\end{equation}
where $\mathcal{L}_{t}$ is the operator associated with the coefficients $A$, $F$ and $G$ in (\ref{stoc}), namely
\begin{equation}\label{L}
\mathcal{L}_{t} \phi(x) = \frac{1}{2}\;{\rm Trace }\left(
G(t,x)G(t,x)^*\nabla^2\phi(x)\right) +
\< x,A^*\nabla\phi(x)\>_H+ \<F(t,x),\nabla\phi(x)\>_H.
\end{equation}

Important results concerning connections between BSDEs and PDEs have been stated by Pardoux and Peng \cite{PaPe2} and  Peng \cite{Pe1}, \cite{Pe2}, \cite{Pe3}, \cite{Pe4} in the Markovian case: these Markovian BSDEs give a Feynman-Kac representation formula for nonlinear parabolic partial differential equations. Conversely, under smoothness assumptions, the solution of a BSDE corresponds to the solution of a system of quasilinear parabolic PDEs.

One of the
main results of this paper, Theorem \ref{main},
specifies conditions
for unique solvability of
equation (\ref{kolmogorovnonlineare-intro}). As usual, if we define the function $v$ by the formula $v(t,x)=Y(t,t,x)$, where $Y(\cdot,t,x)$ is  the solution to BSDE~(\ref{backprimo}), then it turns out that the function
$v$ is deterministic and it is a solution of  the nonlinear
Kolmogorov equation. Various concepts of solutions are known for (linear and) nonlinear
parabolic equations in finite and infinite dimensions. 
Many investigations have been carried out
in connection with the Hamilton Jacobi
Bellman equation arising in optimal control for nonlinear
infinite dimensional stochastic systems.
One
possibility to deal with equation
(\ref{kolmogorovnonlineare-intro}) is to look for classical
solutions, see e.g. \cite{PaPe2} or \cite{Pe2}, i.e. functions which are
twice differentiable with respect to $x$ and once with respect to
$t$, such that $\call [v(t,\cdot)]$ makes sense for every $t\in
[0,T]$ and (\ref{kolmogorovnonlineare-intro}) holds. This requires
heavy assumptions on the functions $\psi$ and $\phi$,
involving existence of derivatives up to order two as well as
trace conditions on second derivatives. Another possibility, in
some sense opposite, is to consider viscosity solutions. Existence
and uniqueness of viscosity solutions   can be proved under much
weaker assumptions on the coefficients in the finite dimensional case, see e.g.
\cite{P2}.  
 However, the theory of viscosity solutions is much harder in the infinite dimensional case.
 Moreover, in view of applications to
optimal control theory, it is important to show the existence of
$\nabla_xv$, since this allows to characterize the optimal control by
feedback laws.  In this paper we will consider  solutions in the so
called mild sense (already considered in the literature, see
\cite{C}, \cite{GoRo} and references within, and, in connection
with the backward stochastic equations approach, \cite{FT1}, \cite{FT2}). Namely a mild
solution $v$ of equation (\ref{kolmogorovnonlineare-intro}) will
satisfy the equality, for $t\in[0,T]$ and $x\in H$,
$$
  v(t,x) =\int_t^TP_{t,\tau}[
\psi (\tau, \cdot,v(\tau,\cdot),G(\tau,\cdot)^*
\nabla_xv(\tau,\cdot)) ](x)\; d\tau+ P_{t,T}[\phi](x),
$$
which arises formally from
 (\ref{kolmogorovnonlineare-intro}) as  the variation of
parameters formula. We  notice that this formula 
is meaningful provided $v$ is only once differentiable with
respect to $x$ and, of course, provided $\psi$, $v$ and
$\nabla_xv$ satisfy appropriate measurability and growth
conditions. Thus, mild solutions are in a sense intermediate
between classical and viscosity solutions. We can prove existence
and uniqueness of a mild solution $v$ by requiring existence and
boundedness (or
 growth conditions) of first derivatives
of $\psi$ and $\phi$: compare Theorem \ref{main}.

As mentioned above, results on system (\ref{backprimo})-(\ref{stoc}) are suitable for applications to problems of nonlinear stochastic optimal control.
Let us consider a controlled Markov process $X^u$ in $H$, on a
time interval $[t,T] \subset [0,T]$, described by an It\^o stochastic
differential equation of the form:
$$
\left \{ \begin{array}{ll}
 dX_s^u= AX_s^u\,ds + F(s,X_s^u)\,ds + G(s,X_s^u)[r(s,X_s^u,u_s)\,ds + dW_s ],& s \in [t,T], \\
    X_t^u=x \in H ,&  \\
\end{array} \right.
$$
where the control process $u$ takes values in a given subset
$\mathcal{U} \subset U$. The aim is to choose a control
process $u$, within a set of admissible controls, in such way to
minimize a cost functional of the form:
$$
J(t,x,u)= \E\left[ \int_t^T \exp \left(\int_t^s
\l(r,X_r^u,u_r)dr \right) l(s,X_s^u,u_s)ds  + \exp \left(\int_t^T \l(r,X_r^u,u_r)dr
\right)\phi(X_T^u)\right],
$$
 where $\l, l, \phi$ are real
functions. $\l$ may be called the discount function.
The Hamilton
Jacobi Bellman equation
 for the value function is of the form (\ref{kolmogorovnonlineare-intro}),
 provided we set
 $$
 \psi(t,x,y,z)= \inf_{u \in \mathcal{U} }\left[
l(t,x,u) + <r(t,x,u),z> + \l (t,x,u)y \right].
$$ 
The control problem is understood in the usual weak sense
 (see \cite{FlSo} and Subsection \ref{ssec-control} below).
 The distinctive feature of our case is the fact that $v$ occurs
explicitly in the Hamilton Jacobi Bellman equation (not only
through its derivatives)  in a nonlinear way
whenever $\l$ effectively depends on $u$ (i.e. we have the control
on the discount). In particular, if we take $r$ bounded and $ \l$
non positive, we obtain a coefficient $\psi$ monotone with respect
to $y$ and Lipschitz with respect to $z$.
Remarks on this special case of Hamilton Jacobi Bellman equation
are in the classical books by Bensoussan \cite{Be} and
Bensoussan-Lions \cite{BeLi}, but they require $\l$ to be bounded.
Fuhrman and Tessitore \cite{FT1}, \cite{FT2} deal with a similar
optimal control problem in infinite dimension, but in their
applications the Hamilton Jacobi Bellman equation depends at
most linearly on $v$.
Under suitable conditions, if we let $v$
denote the unique solution of the Hamilton Jacobi Bellman equation, then we have
$J(t,x,u)\geq v(t,x)$ and the
 equality holds if and only if $u$ and $X^u$ satisfy a
suitable feedback law. Thus, the optimal control $u$ is related to the corresponding optimal trajectory $X^u$
by a feedback law involving $\nabla_x v$.
One of the main results of this paper is the characterization in this
case of the optimal control by a feedback law (Theorem
\ref{maincontrollo}). To this end we use the previous results to guarantee existence, uniqueness and regularity of the solution of Hamilton Jacobi Bellman equation.

The plan of the paper is as follows: Section  \ref{sec-Notation}
is devoted to notations. In Section~\ref{sec-diffbsde}, we state our main result concerning the differentiability of the solution to a BSDE with respect to the data. In Section \ref{secfb}, we apply the previous results  to study the regularity of the map $(t,x)\longmapsto \left(Y(\cdot,t,x),Z(\cdot,t,x)\right)$ solution of the forward-backward system \eqref{backprimo}--\eqref{stoc}. The last section contains the applications to PDEs: in particular, we study the nonlinear Kolmogorov equation and we give some applications to optimal control.

\section{Notations}
\label{sec-Notation}

\subsection{Vector spaces and stochastic processes}

The norm of an element $x$ of a Banach space $E$ will be denoted
$|x|_E$ or simply $|x|$ if no confusion is possible. If
$F$ is another Banach space, $L(E,F)$
denotes the space of bounded linear operators from $E$ to $F$
endowed with the usual operator norm.

The letters $\Xi$, $H$, $K$ will always denote Hilbert spaces.
Scalar product is denoted $\<\cdot,\cdot\>$, with a subscript
to specify the space if necessary. All
Hilbert spaces are assumed to be real and separable. $L_2(\Xi,K)$
is the space of Hilbert-Schmidt operators from $\Xi$ to $K$
endowed with the Hilbert-Schmidt norm.

By a cylindrical Wiener process with values in
  a Hilbert space $\Xi$, defined on a probability
  space $(\Omega, \calf,\P)$, we mean a family
$W(t)$, $t\geq 0$, of linear mappings from $\Xi$ to $L^2(\Omega)$ such that
\begin{description}
  \item[(i)] for every $u\in \Xi$, $\{ W(t)u,\; t\geq 0\}$ is a real
  (continuous) Wiener process;
  \item[(ii)] for every $u,v\in \Xi$ and $t\geq 0$,
  $\E\; (W(t)u\cdot W(t)v)=
  \<u,v\>_\Xi$.
\end{description}

In the following, all stochastic processes will be defined on
subsets of a fixed time interval $[0,T]$.  $\left\{\calf_t\right\}_{t\in [0,T]}$, will denote, except in Subsection \ref{ssec-control}, the natural
  filtration of $W$, augmented with the family $\caln$ of
  $\P$-null sets of $\calf_T$:
$$
\calf_t=\sigma(W(s)\; :\; s\in [0,t])\vee\caln.
$$
The filtration
$\left\{\calf_t\right\}_{t\in[0,T]}$ satisfies the usual conditions. All the concepts of
measurability for stochastic processes (e.g. predictability etc.)
refer to this filtration. For $[a,b]\subset [0,T]$ we also use the notation
$$
\calf_{[a,b]}=\sigma(W(s)-W(a)\; :\; s\in [a,b])\vee\caln.
$$
By $\calp$ we denote the
 predictable $\sigma$-algebra on $\Omega\times [0,T]$
and by $\calb(\Lambda)$ the Borel $\sigma$-algebra of any topological
space $\Lambda$.

Next we define several classes of stochastic processes with values
in a Hilbert space $K$.
\begin{itemize}
  \item
$L_{\calp}^2(\Omega\times [0,T];K)$
denotes the space of equivalence classes of
processes
$Y\in  L^2(\Omega\times [0,T];K)$, admitting a predictable
version. $L_{\calp}^2(\Omega\times [0,T];K)$ is endowed with the norm
$$|Y|^2=\E\left[\int_0^T |Y_\tau|^2\,d\tau\right].
$$

  \item
$\m M^p(K) := L_{\calp}^p(\Omega;L^2([0,T];K))$
denotes the space of equivalence classes of
processes $Z$ such that the norm
$$
\|Z \|_{2,p} =\E\left[\Big(\int_0^T |Z_\tau|^2d\tau\Big)^{p/2}\right]^{1/p}
$$
is finite, and $Z$ admits a predictable
version.

  \item
  $C_{\calp}([0,T];L^2(\Omega;K))$  denotes the space of $K$-valued
processes $Y$ with a predictable
modification such that 
$Y:[0,T]\to L^2(\Omega;K)$ is continuous,
endowed with the norm
$$|Y|^2=\sup_{\tau\in [0,T]}\E\, |Y_\tau|^2.
$$
Elements of $C_{\calp}([0,T];L^2(\Omega;K))$ are
identified up to modification.

\item
$\m S^p(K):=L_{\calp}^p(\Omega; C( [0,T];K))$
  denotes the space of predictable processes
$Y$ with continuous paths in $K$, such that
 the norm
$$
\| Y \|_{\infty,p} =\E\left[\sup\nl_{\tau\in [0,T]}|Y_\tau|^p \right]^{1/p}
$$
is finite.
Elements of $\m S^p(K)$ are
identified up to indistinguishability.

\end{itemize}

If $Y$ is a process in $K$, we will use the following notations:
$$
|Y|_1 = \int_0^T |Y_s|\, ds, \qquad |Y|_\infty = \sup\nl_{t\in[0,T]} |Y_t|, \qquad |Y|_2 = \left(\int_0^T |Y_s|^2 \, ds\right)^{1/2}, \ldots
$$

Given an element $\Psi$ of
$L_{\calp}^2(\Omega\times [0,T];L_2(\Xi,K))$, one can define the It\^o stochastic integral
$\int_0^t\Psi(\sigma)\, d\sigma$, $t\in [0,T]$;
it is a $K$-valued martingale belonging to
 $L_{\calp}^2(\Omega;C( [0,T];K))$.

The previous definitions have obvious extensions to processes
defined on subintervals of $[0,T]$.

\subsection{The class $\calg $}

In this subsection,
$X$, $Y$, $Z$ and $V$ denote Banach spaces.
We recall that
for a mapping $F:X\to V$ the
directional
  derivative at point $x\in X$ in the direction $h\in X$
  is defined as
  $$
  \nabla F(x;h)=\lim_{s\to 0}\frac{F(x+sh)-F(x)}{s},
  $$
  whenever  the limit exists
  in the topology of $V$. $F$ is called G\^ateaux
differentiable at point $x$ if it has directional
derivative in every direction at this point  and
there exists an element of
$L(X,V)$, denoted $\nabla F(x)$ and called G\^ateaux
derivative, such that $\nabla F(x;h)=\nabla F(x)h$
for every $h\in X$.
\begin{definition} 
We say that a mapping $F:X\to V$
belongs to the class $\calg^1 (X;V)$ if it is continuous,
G\^ateaux differentiable on $X$, and  $\nabla F:X\to L(X,V)$
is strongly continuous.
\end{definition}
 The last requirement of the definition means that for every
$h\in X$ the map $\nabla F(\cdot)h:X\to V$ is continuous.
Note that  $\nabla F:X\to L(X,V)$ is not continuous in general
if $L(X,V)$ is
endowed with the norm operator topology; clearly,
if this happens
then $F$ is
Fr\'echet differentiable on $X$.
Some features of the class $\calg^1(X,V)$ are collected below.

\begin{lemma}\label{proprietadig}
  Suppose $F\in\calg^1(X,V)$. Then
\begin{description}
  \item[(i)] $(x, h)\mapsto \nabla F(x)h$ is
  continuous from $X\times X$ to $V$;
  \item[(ii)]  If $G\in
  \calg^1(V,Z)$ then $G( F)\in \calg^{1} (X,Z)$ and
  $\nabla(G( F))(x)=\nabla G(F(x))\nabla F(x)$.
\end{description}
\end{lemma}

\begin{lemma}
\label{staing}
  A map $F:X\to V$ belongs to $\calg^1(X,V)$
  provided the following conditions hold:
\begin{description}
  \item[(i)]  the directional
  derivatives $\nabla F(x;h)$ exist
  at every point $x\in X$ and in every direction $h\in X$;
  \item[(ii)] for every  $h$, the mapping
   $\nabla F(\cdot;h): X\to V$
  is continuous;
  \item[(iii)] for every $x$, the mapping $h\mapsto\nabla F(x;h)$
  is continuous from $X$ to $V$.
\end{description}
\end{lemma}

The proofs of these lemmas are in \cite{FT1}. We need to generalize these definitions to functions depending on
several variables.  For a function $F:X\times Y\to V$ the partial
directional and G\^ateaux derivatives
with respect to the first argument, at point $(x,y)$ and in the
direction  $h\in X$, are
denoted
$\nabla_x F(x,y;h)$ and $\nabla_x F(x,y)$ respectively,
their definitions being obvious.

\begin{definition}
We say that a mapping $F:X\times Y\to V$
belongs to the class $\calg^{1,0} (X\times Y;V)$ if it is continuous,
G\^ateaux differentiable with respect to
$x$ on $X\times Y$, and  $\nabla_xF:X\times Y\to L(X,V)$
is strongly continuous.
\end{definition}

As in Lemma \ref{proprietadig} 
for  $F\in\calg^{1,0}(X\times Y,V)$ the map
$(x, y, h)\mapsto \nabla_xF(x,y)h$ is
  continuous from $X\times Y\times X$ to $V$, and
analogues of the previously stated  chain rules hold.
The following result is proved as Lemma
\ref{staing} (but
note that continuity is explicitly required).

\begin{lemma}\label{staing2}
  A continuous map $F:X\times Y\to V$ belongs to
  $\calg^{1,0}(X\times Y,V)$
  provided the following conditions hold:
\begin{description}
  \item[(i)]
  the directional
  derivatives $\nabla_x F(x,y;h)$ exist
  at every point $(x,y)\in X\times Y$
  and in every direction $h\in X$;
  \item[(ii)] for every  $h$, the mapping
   $\nabla F(\cdot,\cdot ;h): X\times Y\to V$
  is continuous;
  \item[(iii)] for every $(x,y)$,
   the mapping $h\mapsto\nabla_x F(x,y;h)$
  is continuous from $X$ to $V$.
\end{description}
\end{lemma}

When $F$ depends on additional arguments,
 the previous definitions and properties  have
obvious generalizations. For instance, we say that
$F:X\times Y\times Z\to V$
belongs to  $\calg^{1,1,0} (X\times Y\times Z;V)$
if it is continuous,
G\^ateaux differentiable with respect to
$x$ and $y$ on $X\times Y\times Z$,
and  $\nabla_xF:X\times Y\times Z\to L(X,V)$ and
$\nabla_yF:X\times Y\times Z\to L(Y,V)$ are strongly continuous.

\section{Differentiability of BSDEs}
\label{sec-diffbsde}

In this section we want to study the BSDE
\begin{equation}
\label{beqp}
 Y_t =\xi +\int_t^T \psi(s,X_s,Y_s,Z_s)\,ds - \int_t^T Z_s\,dW_s,\qquad 0\leq t\leq T
\end{equation}
where $\xi$ is a given $\mathcal{F}_T$--measurable random variable with values in $K$ and $X_{\t}$, $\t \in [0,T]$, is a given continuous predictable process with values in $H$. To stress the dependence of the solution to the BSDE with respect to the terminal condition $\xi$ and the process $X$, we will denote $\left(Y(X,\xi), Z(X,\xi)\right)$ the solution to \eqref{beqp}. The main point is to study the regularity of the map $(X,\xi)\longmapsto \left(Y(X,\xi), Z(X,\xi)\right)$. 

$\psi: [0,T]\times H \times K \times L_2(\Xi,K) \rightarrow K$ is assumed to be a Borel--measurable function which satisfies moreover the following assumption:

\begin{hypothesis}\label{ipsupsi}
There exist some constants $m\geq 0$, $n\geq 0$, $c\geq 0$, $L\geq 0$ and $\mu\in\rset$ such that
\begin{description}
 \item[(i)] for every $\sigma\in [0,T]$, $x\in H$, $y \in K$, $z_1,z_2\in L_2(\Xi,K)$,
    $$
    |\psi(\sigma,x, y,z_1)-\psi(\sigma,x, y,z_2)| \leq L\,|z_1-z_2|;
  $$
  \item[(ii)]  for every $\sigma\in [0,T]$, $x\in H$,  $z\in L_2(\Xi,K)$, $y\longmapsto \psi(\sigma,x,y,z)$ is $\mu$--monotone meaning that
$$
\forall (y_1,y_2)\in K^2,\qquad \<y_1-y_2, \psi(\sigma,x,y_1,z)- \psi(\sigma,x,y_2,z)\> \leq \mu |y_1-y_2|^2~;
$$
 \item[(iii)] for every $\sigma\in [0,T]$,
  $\psi (\sigma,\cdot,\cdot,\cdot)\in \m G^{1,1,1} (H \times K \times
 L_2(\Xi,K))$.
  \item[(iv)] for every $\sigma\in [0,T]$, $x \in H$, $
  y\in K$, $z\in L_2(\Xi,K)$,
  $$
  |\nabla_x\psi(\sigma,x, y,z)|+|\nabla_y\psi(\sigma,x, y,z)|\leq
  q(\sigma)+c\left(|x|^m+|y|^n+ |z|^2\right)
  $$
  where $q\in \lp^1(0,T)$.
\item[(v)] $\psi(s):=\psi(s,0,0,0)$ belongs to $L^1(0 ,T)$.
\end{description}
\end{hypothesis}

\begin{remark} It follows from Hypothesis \ref{ipsupsi} i), iii) and iv) that 
\begin{eqnarray*}
\lefteqn{\left|\psi(s,x_1,y_1,z_1)-\psi(s,x_2,y_2,z_2)\right| }
\\
&& \leq   L\,|z_1-z_2| + \left(|x_1-x_2|+|y_1-y_2| \right) \left( q(s) + C\left(|x_1|^m+|x_2|^m + |y_1|^n + |y_2|^n + |z_1|^2\right)\right) \nonumber.
\end{eqnarray*}

Consequently, we have
\begin{eqnarray}
\label{majpsi}
\left|\psi(s,x,y,z)\right| \leq L|z| +(|x| +|y|)(q(s) +C( |x|^m +|y|^n )) + |\psi(s,0,0,0)|.
\end{eqnarray}
\end{remark}

For $p \geq 1$ we denote $\mathcal{K}^p$ the Banach space
$$
\mathcal{K}^p = \m S^p(K) \times \mathcal{M}^p\left( L_2(\Xi,K) \right)
$$
endowed with the  norm
$$
\left\| ( Y,Z ) \right\|_p= \left( \left\| Y \right\|_{\infty,p}^p +\left\|Z\right\|_{2,p}^p\right)^{1/p}= \E\left[ \sup\nl_{t \in [0,T]} |Y_t|^p + \Big(\int_0^T ||Z_t||^2\, dt \Big)^{p/2}\right]^{1/p}.
$$

\begin{proposition}
\label{solbak} 
Let the assumption \ref{ipsupsi} hold.
 
Let $p>1$ and let $\xi\in\lp^p(\Omega;K)$, $X\in\m S^{p(m+1)}(H)$. The BSDE \eqref{beqp} has a unique solution in $\calk^p$, $(Y(X,\xi),Z(X,\xi))$, which satisfies
\begin{equation}
   ||\left( Y(X,\xi),Z(X,\xi)\right)||_p \leq
C\left(1 + ||\xi||_{p}+ ||X||_{\infty,p(m+1)}^{ m+1}\right)
 \label{stimaY}
\end{equation}
for a suitable constant $C$ depending only on $p$, $T$, $L$ and $\mu$.
\end{proposition}

\begin{proof}
The proof relies heavily on the results given in \cite{BDHPS} in the finite dimensional case. The generalization of these results to the case of a cylindrical Wiener process taking its values in an Hilbert space is immediate.

\medskip

Let $X \in\m S^{p(m+1)}(H)$ and $\xi\in\lp^p(\Omega;K)$. Since the function $\psi$ is Lipschitz with respect to $z$ and $\mu$--monotone with respect to $y$, Theorem 4.2 in \cite{BDHPS} shows that the previous BSDE has a unique solution, $(Y,Z)$ in $\m K^p$. Moreover, it follows from \cite{BDHPS}[Proposition 3.2]  that
$$
\E\left[ \sup\nl_{\tau\in[0,T]}|Y_\tau|^p+ \left(\int_0^T|Z_\sigma|^2d
\sigma\right)^{p/2} \right]\leq K\, \E \left[ |\xi|^p +\left(\int_0^T \vert \psi(\sigma,X_{\s},0,0) \vert d
\sigma\right)^p \right].
$$
Hence, taking into account \eqref{majpsi} we have, 
$$
\left\| (Y,Z) \right\|_p \leq K \left( 1 + ||\xi||_p + \| X\|_{\infty, (m+1)p}^{m+1} \right)
$$
which is exactly \eqref{stimaY}.
\end{proof}

Let us introduce some further notations. For any $p>1$, we denote by $\dsp$ the product space $\m S^{p(m+1)}(H) \times L^p(\Omega;K)$ and we consider, according the previous proposition, the map $\Phi$ from $\dsp$ to $\m K^p$ defined by $\Phi(X,\xi)=\left(Y(X,\xi),Z(X,\xi) \right)$ where $\left(Y(X,\xi),Z(X,\xi)\right)$ stands for the solution to the BSDE~\eqref{beqp}.

\begin{proposition}
\label{diffprop}
Let the assumption \ref{ipsupsi} hold and let $p_*= \max(2,n) + \frac{1}{m+1}$.

\begin{description}
\item[(a)] If $p > p_*$, the map $\Phi$ is continuous from $
\dsp$ to $\mathcal{K}^{r}$ where  $r = p/p_*$.

\item[(b)] If $p > p_* + \max(np_*,2)$,  $\Phi$ belongs to $\mathcal{G}^{1}(\dsp,\mathcal{K}^{\rho})$ with $r=p/p_*$ and $\rho= r \min\left(\frac{1}{n+1},\frac{p_*}{2+p_*}\right)$.

Moreover, for all $(X,\xi) \in \dsp$,
 the directional derivative of $\Phi$ in the direction $(N,\zeta)\in \dsp$, $\nabla_{X,\xi}\Phi(X,\xi)(N,\zeta)$ is the unique solution $(G,H)$ in $\mathcal{K}^{r}$ to :
\begin{eqnarray}
\nonumber
G_{t}+\int_{t}^{T} H_{\sigma}dW_{\sigma} & = & \zeta + \int_{t}^{T}\nabla_{x}
\psi(\sigma,X_{\sigma},Y_{\sigma}(X,\xi),Z_{\sigma}(X,\xi))N_{\sigma }\,d
\sigma     \\
\label{bsdegradient}
&& + \int_{\tau}^{T}\nabla_{y}\psi(\sigma, X_{\sigma},Y_{\sigma}(X,\xi),Z_{\sigma}(X,\xi))
G_\sigma\, d \sigma \\
\nonumber
&& + \int_{\tau}^{T}\nabla_{z}
\psi(\sigma,X_{\sigma},Y_{\sigma}(X,\xi),Z_{\sigma}(X,\xi))
H_\sigma\,d \sigma.
\end{eqnarray}

\item Finally the following estimate holds:
\begin{equation}
\label{stimanablaYZ}
     \left\| \nabla_{X,\xi}\Phi(X,\xi)(N,\zeta)
\right\|_r \leq K \left( \|\zeta\|_p + \| N \|_{\infty,p(m+1)} \left(1 + \|\xi\|_p^{n\vee 2} + \| X \|_{\infty,p(m+1)}^{(m+1)(n\vee 2)} \right) \right).
\end{equation}
\end{description}
\end{proposition}

\begin{proof}

\textbf{(a)} Let us assume that $p > p_*$ and let us set $r= p/p_* >1$. Using Proposition 3.2 in \cite{BDHPS}, we have, if $(Y,Z)=\Phi(X,\xi)$ and $(Y',Z')=\Phi(X',\xi')$ where $(X,\xi)$ and $(X',\xi')$ belongs to $\dsp$, 
\begin{eqnarray*}
\left\| (Y,Z) - (Y',Z') \right\|_r^r & \leq & C\, \E\left[ |\xi-\xi'|^r + \left(\int_0^T \left|\psi(s,X_s,Y_s,Z_s)-\psi(s,X'_s,Y_s,Z_s)\right| ds\right)^r \right] \\
& \leq &   C\, \E\left[ |\xi-\xi'|^r + \left| X-X' \right|_\infty^r \left( |q|_1^r + \left|X \right|_\infty^{mr} + \left|X'\right|_\infty^{mr} + \left| Y \right|_\infty^{nr} + \left| Z \right|_2^{2r} \right) \right] \\
& \leq &   C\, \E\left[ |\xi-\xi'|^r + \left| X-X' \right|_\infty^r \left( 1 + \left|X \right|_\infty^{mr} + \left|X'\right|_\infty^{mr} + \left| Y \right|_\infty^{nr} + \left| Z \right|_2^{2r} \right) \right].
\end{eqnarray*}
Using H\"older's inequality, we deduce that
%
%
$$
\left\| (Y,Z) - (Y',Z') \right\|_r  \leq  C\left( || \xi-\xi'||_r  + \left\| X-X' \right\|_{\infty,p(m+1)} \times A \right)
$$
with $A$ given by
$$
A = 1 + \| X \|_{\infty, u(m,r)}^m + \left\| X'\right\|_{\infty, u(m,r)}^m + \| Y \|_{\infty, u(n,r)}^n + \| Z\|_{2,u(2,r)}^2
$$
where, for any $a>0$, $u(a,r)= a \frac{rp(m+1)}{p(m+1)-r} = a \frac{p}{2\vee n}$.  It follows that $\left\|(Y,Z) - \left(Y',Z'\right)\right\|_r$ is bounded from above by, up to a multiplicative constant $C$,
$$
\left\| \xi-\xi'\right\|_p + \left\| X-X'\right\|_{\infty, p(m+1)} \left(1 + \| X \|_{\infty, p(m+1)}^m + \left\| X' \right\|_{\infty, p(m+1)}^m + \| Y \|_{\infty, p}^n + \|Z\|_{2, p}^2\right) .
$$
This gives the continuity of $\Phi$.

\medskip

\textbf{(b)} Let us assume now that $p > p_* + \max (np_*,2)$ and let us define  $r= p/p_*$ together with $\rho= r \min\left(\frac{1}{n+1},\frac{p_*}{2+p_*}\right) = \min\left(\frac{p}{p_*(n+1)}, \frac{p}{2+p_*} \right) > 1$. 

Let us pick $(X,\xi)$ and $(N,\zeta)$ in $\dsp$ and, for simplicity, let us denote by $T$ the triple $(X,Y,Z)=(X,\Phi(X,\xi))$. We consider the BSDE
$$
G_t  =  \zeta  + \int_t^T \nabla_x\psi(s,T_s) N_s \, ds + \int_t^T \nabla_y\psi(s,T_s) G_s \, ds + \int_t^T \nabla_z\psi(s,T_s) H_s \, ds - \int_t^T H_s\,dB_s .
$$
Let us observe that
$$
\int_0^T \left|  \nabla_x\psi(s,T_s) N_s  \right| \, ds\leq K |N|_\infty \left( |q|_1 + |X|_\infty^m + |Y|_\infty^n + |Z|_2^2 \right)
$$
and arguing as before with H\"{o}lder inequality we deduce that
$$
\left[\E \left(\int_0^T \left| \nabla_x\psi(s,T_s) N_s\right| \, ds\right) ^r \right]^{1/r} \leq K \| N \|_{\infty,p(m+1)} \left(1 + \| X \|_{\infty,p(m+1)}^m + \| Y \|_{\infty,p}^n + \| Z\|_{2,p}^2 \right).
$$
Moreover, we have, since $\psi$ is Lipschitz with respect to $z$ and $\mu$--monotone with respect to $y$,
$$
\left| \nabla_z\psi(s,T_s) \left(h-h'\right) \right| \leq L \left|h-h'\right|,\quad \left( g-g' , \nabla_y\psi(s,T_s) \left(g-g'\right) \right)_K \leq \mu \left|g-g'\right|^2.
$$
It follows from \cite{BDHPS} that the previous BSDE has a unique solution $(G,H)$ in $\m K^r$ such that
$$
\|(G,H)\|_r \leq K \left( \|\zeta\|_r + \| N \|_{\infty,p(m+1)} \left(1 + \| X \|_{\infty,p(m+1)}^m + \| Y \|_{\infty,p}^n + \| Z\|_{2,p}^2 \right) \right).
$$
Taking into account the inequality \eqref{stimaY}, we deduce that
\begin{equation}
\label{sc}
\|(G,H)\|_r \leq K \left( \|\zeta\|_p + \| N \|_{\infty,p(m+1)} \left(1 + \|\xi\|_p^{n\vee 2} + \| X \|_{\infty,p(m+1)}^{(m+1)(n\vee 2)} \right) \right).
\end{equation}

\smallskip

It remains to prove that the directional derivative of $\Phi$ at $(X,\xi)\in\dsp$ in the direction $(N,\zeta)\in\dsp$ is given by $(G,H)$. For $\ep>0$,  we set $X^\ep = X + \ep N$ and we consider $(Y^\ep,Z^\ep)$ the solution in $\m K^p$ to the BSDE
$$
Y^\ep_t = \xi +\ep\,\zeta + \int_t^T \psi(s,X^\ep_s,Y^\ep_s,Z^\ep_s)\,ds - \int_t^T Z^\ep_s\,dB_s.
$$
Before going further, let us introduce the following notations: $U^\ep_t = \ep^{-1}\left(Y^\ep_t - Y_t\right) - G_t$ and $V^\ep_t = \ep^{-1}\left(Z^\ep_t - Z_t\right) - H_t$. We have,
\begin{eqnarray*}
U^\ep_t  & = &  \frac{1}{\ep} \int_t^T \left(\psi(s,X^\ep_s,Y^\ep_s,Z^\ep_s) - \psi(s,X_s,Y_s,Z_s)\right) ds  - \int_t^T V^\ep_s\,dB_s  \\
& & - \int_t^T \nabla_x\psi(s,T_s) N_s \, ds - \int_t^T \nabla_y\psi(s,T_s) G_s \, ds - \int_t^T \nabla_z\psi(s,T_s) H_s \, ds.
\end{eqnarray*}
Using the fact that $\psi(s,\cdot,\cdot,\cdot)$ belongs to $\m G^{1,1,1}$, we can write
$$
\frac{1}{\ep}\left(\psi(s,X^\ep_s,Y^\ep_s,Z^\ep_s) - \psi(s,X_s,Y_s,Z_s)\right) = A^\ep_{s} \, \frac{Y^\ep_{s}-Y_s}{\ep} + B^\ep_s\, \frac{Z^\ep_s - Z_s}{\ep} + C^\ep_s
$$
where $A^\ep_s\in L(K,K)$ and $B^\ep_s\in L\left( L_2(\Xi,K), K\right)$ are defined by
$$
\forall y\in K,\qquad A^\ep_s y = \int_0^1 \nabla_y\psi\left(s,X^\ep_s,Y_s +\alpha (Y^\ep_s-Y_s), Z_s\right) y \,d\alpha,
$$
$$
\forall z\in L_2(\Xi,K),\qquad B^\ep_s z = \int_0^1 \nabla_z\psi\left(s,X^\ep_s,Y^\ep_s,Z_s +\alpha (Z^\ep_s-Z_s)\right) z \,d\alpha
$$
and where 
$$
C^\ep_s = \frac{1}{\ep}\left(\psi(s,X^\ep_s,Y_s,Z_s) - \psi(s,X_s,Y_s,Z_s)\right) = \int_0^1 \nabla_x \psi(x,X_s + \alpha\ep N_s, Y_s,Z_s)N_s\,d\alpha.
$$
Then $(U^\ep,V^\ep)$ solves the following BSDE
$$
U^\ep_t = \int_t^T \left( A^\ep_s U^\ep_s + B^\ep_s V^\ep_s\right) ds + \int_t^T \left( P^\ep(s) + Q^\ep(s) + R^\ep(s) \right) ds - \int_t^T V^\ep_s\,dB_s
$$
where we have set
$$
P^\ep(s) = \left( A^\ep_s - \nabla_y \psi(s,T_s) \right) G_s, \quad Q^\ep(s) = \left( B^\ep_s - \nabla_z \psi(s,T_s) \right) H_s,
$$
together with
$$
R^\ep(s)  =  \int_0^1 \nabla_x \psi(x,X_s + \alpha\ep N_s, Y_s,Z_s)N_s\,d\alpha -\nabla_x \psi(s,T_s) N_s.
$$

Since $\left(y,A^\ep_s\,y\right) \leq \mu |y|^2$ and $| B^\ep_s z |\leq L |z|$, we can apply Proposition 3.2 in \cite{BDHPS} to learn that, 
$$
\left\| (U^\ep,V^\ep) \right\|_\rho \leq K \left( \left\| |P^\ep|_1 \right\|_\rho + \left\| |Q^\ep|_1 \right\|_\rho + \left\| |R^\ep|_1 \right\|_\rho\right).
$$
for a suitable constant $K$ depending on $L$, $\mu$, $T$ and $\rho$ and where we use the notation 
$$
\left| P^\ep \right|_1 = \left| P^\ep(\cdot) \right|_1 = \int_0^T \left| P^\ep(s) \right| ds.
$$
It remains to check that the right hand side goes to zero as $\ep$ tends to 0.

Let us recall that
\begin{eqnarray*}
P^\ep(s)&  = & \int_0^1 \nabla_y \psi(s,X^\ep_s,Y_s+\alpha(Y^\ep_s-Y_s), Z_s  )G_s\, d\alpha - \nabla_y \psi(s,T_s) G_s \\
& = & \int_0^1 \left( \nabla_y \psi(s,X^\ep_s,Y_s+\alpha(Y^\ep_s-Y_s), Z_s  ) - \nabla_y \psi(s,T_s)\right) G_s \, d\alpha, 
\end{eqnarray*}
so that
$$
\left| P^\ep(s) \right| \leq \int_0^1 \left| u^\ep(\alpha,s) \right| d\alpha
$$
with $u^\ep(\alpha,s)=\left( \nabla_y \psi(s,X^\ep_s,Y_s +\alpha(Y^\ep_s-Y_s), Z_s ) - \nabla_y \psi(s,T_s) \right) G_s$.

We consider the product space $[0,1]\times[0,T]\times \Omega$ with the measure $\lambda\otimes\lambda\otimes\P$ where $\lambda$ is the Lebesgue measure. Since $\psi(s,\cdot,\cdot,\cdot)$ belongs to  $\m G^{1,1,1}$, $u^\ep(s,\alpha)$ goes to zero in $\lambda\otimes\lambda\otimes\P$--measure by continuity of $\Phi$ and we have,
\begin{eqnarray*}
\left| u^\ep(\alpha,s) \right| & \leq & K\left( q(s) + |Z_s|^2 + |X_s|^m + |X_s^\ep|^m + |Y_s|^n + |Y_s^\ep|^nÊ\right) |G_s| \\
& \leq & K\left( q(s) + |Z_s|^2 + |X|_\infty^m + |X^\ep|_\infty^m + |Y|_\infty^n + |Y^\ep|_\infty^nÊ\right) |G|_\infty .
\end{eqnarray*}
First of all, let us recall that $|G|_\infty$ belongs to $\lp^r(\Omega)$ and that $q\in\lp^1(0,T)$ and $|Z|_2^2$ belongs to $\lp^{p/2}$.
Moreover, $Y^\ep$ converges to $Y$ in $\m S^r(K)$, $X^\ep$ converges to $X$ in $\m S^{p(m+1)}(H)$ so that  $\left\{ \left|Y^\ep \right|_\infty^n \right\}_{\ep>0}$ is uniformly integrable in $\lp^{r/n}$ and $\left\{ \left| X^\ep\right|_\infty^m \right\}_{\ep >0}$ is uniformly in $\lp^{p(m+1)/m}$. Since we have $\frac{1}{r} + \frac{2}{p} = \frac{2+p_*}{p} \leq 1$, $\frac{n}{r}+\frac{1}{r} = \frac{(n+1)p_*}{p} \leq 1$ and $\frac{1}{r} + \frac{m}{p(m+1)} = \frac{\max(2,n) + 1}{p} \leq 1$, the right hand side of the previous inequality is uniformly integrable as a function of the three variables $(\alpha,s,\omega)$. In particular,  $\left| P^\ep \right|_1$ converges to 0 in probability as $\ep$ goes to zero. Moreover we have
$$
\left| P^\ep \right|_1 \leq K\left(1+ |Z|_2^2 + |X|_\infty^m + |X^\ep|_\infty^m + |Y|_\infty^n + |Y^\ep|_\infty^n \right) |G|_\infty.
$$
Since $|G|_\infty$ is in $\lp^r$ we have only to check that the RHS of the previous inequality is uniformly integrable in $\lp^{\rho\sigma}$ where $\sigma = r/(r-\rho)$ is the conjugate exponent of $r/\rho$ to prove that $\left| P^\ep \right|_1$ tends to 0 in $\lp^\rho$ as $\ep$ goes to 0. Let us observe that 
$$
\rho\sigma = r\,\frac{\min\left(\frac{1}{n+1},\frac{p_*}{2+p_*}\right)}{1-\min\left(\frac{1}{n+1},\frac{p_*}{2+p_*}\right)} = r \, \min\left(\frac{1}{n},\frac{p_*}{2} \right) = \min\left( \frac{p}{np_*}, \frac{p}{2}\right).
$$
As we said before $|Z|_2^2\in\lp^{p/2} \subset \lp^{\rho\sigma}$, $|Y^\ep|_\infty^n$ is uniformly integrable in $\lp^{r/n}$ and $ r/n = {p/np^*}\geq \rho\sigma$ and $|X^\ep|_\infty^m$ is uniformly integrable in $\lp^{p(m+1)/m}$ and $\rho\sigma \leq p/2 \leq p(1+1/m)$. Thus $\left| P^\ep \right|_1$ goes to 0 in $\lp^\rho$.

\smallskip

The term $R^\ep$ can be treated exactly in the same way as we did for $P^\ep$. Indeed, we have
\begin{eqnarray*}
R^\ep(s)&  = & \int_0^1 \nabla_x \psi(s,X_s+\alpha\ep N_s,Y_s, Z_s  )N_s \, d\alpha  - \nabla_x \psi(s,T_s) N_s \\
& = & \int_0^1 \left( \nabla_x \psi(s,X_s+\alpha\ep N_s,Y_s, Z_s  ) - \nabla_x \psi(s,T_s)\right)    N_s\, d\alpha,
\end{eqnarray*}
$\nabla_x \psi$ satisfies the same growth condition as $\nabla_y \psi$ and  $N$ belongs to $\m S^{p(m+1)}(H)$ which is a better situation than the previous one where we have only $G\in\m S^r(K)$.

\smallskip

Let us see that $\left|Q^\ep\right|_1$ goes to 0 in $\lp^\rho$ as $\ep$ goes to 0. We write
\begin{eqnarray*}
Q^\ep(s)&  = & \int_0^1 \nabla_z \psi(s,X^\ep_s,Y^\ep_s, Z_s +\alpha (Z^\ep_s-Z_s) ) H_s\, d\alpha  - \nabla_z \psi(s,T_s) H_s \\
& = & \int_0^1 \left( \nabla_z \psi(s,X^\ep_s,Y^\ep_s, Z_s +\alpha (Z^\ep_s-Z_s) ) - \nabla_z \psi(s,T_s)\right) H_s\, d\alpha.
\end{eqnarray*}
Since $\psi(s,\cdot,\cdot,\cdot)$ is $\m G^{1,1,1}$, $v^\ep(s,\alpha):= \left( \nabla_z \psi(s,X^\ep_s,Y^\ep_s, Z_s +\alpha (Z^\ep_s-Z_s) ) - \nabla_z \psi(s,T_s)\right) H_s$ goes to zero in $\lambda\otimes\lambda\otimes\P$--measure and since $\psi$ is Lipschitz with respect to $z$, we have $\left| v^\ep(s,\alpha) \right | \leq 2L \, |H_s|$. Taking into account the fact that $H$ belongs to the space $\mathcal{M}^r\left( L_2(\Xi,K) \right)$, we deduce that $v^\ep$ goes to 0 in $\lp^1((0,1)\times(0,T)\times\Omega)$. In particular, $\left|Q^\ep\right|_1$ goes to 0 in probability. 
We have also $\left|Q^\ep\right|_1 \leq 2L\sqrt{T}\, |H|_2$ and since $|H|_2$ belongs to $\lp^r\subset \lp^\rho$, the bounded convergence theorem gives $\lim_{\ep\to 0} \left\| \left|Q^\ep \right|_1 \right\|_\rho = 0$.

Thus we have proved that $(G,H) = \nabla_{(X,\xi)} \Phi(X,\xi)(N,\zeta)$ in $\m K^\rho$.

\medskip

By linearity, it follows directly from \eqref{sc}, that $(N,\zeta)\fl (G,H)$ is continuous from $\dsp$ to $\m K^r$. Let us prove that $(X,\xi)\fl (G,H)$ is continuous from $\dsp$ to $\m K^\rho$.  Let us consider $(G,H)\in\m K^r$ and $(G',H')\in\m K^r$ the solutions to \eqref{bsdegradient} associated to $(X,\xi)\in\dsp$ and $(X',\xi')\in\dsp$ ($(N,\zeta)$ is fixed in $\dsp$). Once again, using the a priori estimates in \cite{BDHPS}, we have 
\begin{equation}
\label{ouf}
\left\|(G,H)-(G',H') \right\|_\rho \leq C  \left( \left\| |P|_1 \right\|_\rho + \left\| |Q|_1 \right\|_\rho + \left\| |R|_1 \right\|_\rho\right)
\end{equation}
for a suitable constant $C$  where we have set
\begin{eqnarray*}
P(s)  & =  & \left\{ \nabla_y \psi\left(s,X'_s,Y'_s, Z'_s  \right) - \nabla_y \psi(s,X_s,Y_s, Z_s  ) \right\} G_s, \\
Q(s)  & =  & \left\{ \nabla_z \psi\left(s,X'_s,Y'_s, Z'_s  \right) - \nabla_z \psi(s,X_s,Y_s, Z_s  ) \right\} H_s, \\
R(s)  & =  & \left\{ \nabla_x \psi\left(s,X'_s,Y'_s, Z'_s  \right) - \nabla_x \psi(s,X_s,Y_s, Z_s  ) \right\} N_s.
\end{eqnarray*}
If $\left(X',\xi'\right)\fl (X,\xi)$ in $\dsp$ then $\left(Y',Z' \right)\fl (Y,Z)$ in $\m K^r$ and, since $\psi(s,\cdot,\cdot,\cdot)$ belongs to $\m G^{1,1,1}$, arguing exactly as we did to show the G\^ateau-differentiability of $\Phi$, the right hand side of \eqref{ouf} tends to 0.

\medskip

In conclusion, $\Phi$ belongs to $\m G^1\left( \dsp, \m K^\rho \right)$.
\end{proof}

\begin{remark}
This result can be rewritten in the following way.

Let $\rho >1$ and let us define $p = \rho\,\max\left( (n+1) p_*, p_* +2 \right) $ where $p_* = \max(2,n) + \frac{1}{m+1}$. Then the map  $(X,\xi)\fl (Y(X,\xi),Z(X,\xi))$ from $\dsp[p]$ to $\m K^\rho$ is in $\m G^1$.
\end{remark}

\section{The forward-backward system}
\label{secfb}

In this section, we apply the previous results on the differentiability of BSDEs to study the differentiability of the solution to the forward-backward system \eqref{backprimo}--\eqref{stoc}.

We start by recalling some results on the solution to the forward equation.

\subsection{The forward equation}
\label{sec-forward}

Let us recall that $\left\{ W_t \right\}_{t\in [0,T]}$ is  a cylindrical Wiener process with values in a Hilbert space $\Xi$, defined on a probability space $(\Omega, \calf,\P)$.
We fix an interval $[t,T]\subset [0,T]$ and
we consider the It\^o stochastic
differential equation:
\begin{equation}
\label{equazioneforward}
\left\{\begin{array}{l}\dis
dX_\tau = AX_\tau\; d\tau+F(\tau,X_\tau)\; d\tau
+G(\tau,X_\tau)\; dW_\tau,\qquad \tau\in [t,T],\\
\dis X_t=x\in H.
\end{array}\right.
\end{equation}

We assume the following:
\begin{hypothesis}
\label{ipotesiuno}
\begin{description}
  \item[(i)] The operator $A$ is the generator of a strongly
  continuous semigroup $e^{tA}$, $t\geq 0$, in the Hilbert space
  $H$.
  \item[(ii)] The mapping $F:[0,T]\times H\to H$ is measurable
  and satisfies, for
  some constant $L>0$,
 \begin{eqnarray*}
  |F(t,x)-F(t,y)| & \leq & L\, |x-y|,\qquad t\in [0,T],\; x,y\in H, \\
|F(t,x)| & \leq & L\,   (1+|x|), \qquad t\in [0,T],\; x \in H.
\end{eqnarray*}
  \item[(iii)] $G: [0,T]\times H \fl L(\Xi,H)$ is such that,
  for every $v\in \Xi$, the map $Gv: [0,T]\times H\to H$ is measurable,
  $e^{sA}G(t,x)\in L_2(\Xi,H)$ for every $s>0$,
  $t\in [0,T]$ and $x\in H$,
   and
\begin{eqnarray*}
|e^{sA}G(t,x)|_{L_2(\Xi,H)} & \leq & L\; s^{-\gamma}
  (1+|x|),\\
  |e^{sA}G(t,x)-e^{sA}G(t,y)|_{L_2(\Xi,H)} & \leq  & L\; s^{-\gamma}
  |x-y|, \\
|G(t,x)|_{L(\Xi,H)} & \leq & L\;   (1+|x|), 
\end{eqnarray*}
 for
  some constants $L>0$ and $\gamma\in [0,1/2)$.
\item[(iv)] For every $s>0$,
$t\in [0,T]$,
$$F(t,\cdot)
\in \calg^1(H, H),\qquad
e^{sA}G(t,\cdot)\in\calg^1(H, L_2(\Xi,H)).
$$
  \end{description}
\end{hypothesis}

By a solution of equation (\ref{equazioneforward}) we mean
an $(\calf_t)$-predictable process $X_\tau$, $\tau\in [t,T]$,
with continuous paths in $H$, such
that, $\P$-a.s.
\begin{equation}
\label{equazioneforwardmild}
X_\tau = e^{(\tau-t)A}x+\int_t^\tau e^{(\tau-\sigma)A}
F(\sigma,X_\sigma)\; d\sigma
+\int_t^\tau e^{(\tau-\sigma)A}
G(\sigma,X_\sigma)\; dW_\sigma,\qquad \tau\in [t,T].
\end{equation}

To stress the dependence on initial data, we denote the solution by $X(\tau,t,x)$.
Note that $X(\tau,t,x)$ is $\calf_{[t,T]}$-measurable, hence independent of $\calf_t$.
A consequence of the previous assumptions is that, for every $s>0$, $t\in [0,T]$, $x,h\in H$,
 $$
  |\nabla_x F(t,x)h|\leq L\; |h|, \qquad
|\nabla_x(e^{sA}G(t,x))h|_{L_2(\Xi,H)}\leq L\; s^{-\gamma} |h|.
$$

The following results are proved by Fuhrman and Tessitore in \cite{FT1}.
\begin{proposition}
\label{esisteunicox}
Under the assumptions of Hypothesis
 \ref{ipotesiuno}-$(i)$-$(ii)$-$(iii)$, \eqref{equazioneforward} has a unique solution $X$. 
For every $p\geq 1$, $X$ belongs to $\m S^p(H)$ and 
\begin{equation}
\label{stimadeimomentidix}
  \E\left[ \sup\nl _{\tau\in[t,T]}|X_\tau |^p \right]\leq C(1+|x|)^p,
\end{equation}
for some constant $C$ depending only on $p,\gamma,T,L$ and
$M:=\sup_{\tau\in [0,T]}|e^{\tau A}|$.
\end{proposition}

To go further, we need to investigate the dependence of the solution
$X(\tau,t,x)$ on the initial data $x$ and $t$.
We first
reformulate equation \eqref{equazioneforwardmild}  as an
equation on $[0,T]$. We set
$$
  S(\tau)=e^{\tau A} {\rm \quad
for \quad} \tau\geq 0, \qquad\qquad S(\tau)=I
{\rm \quad
 for\quad} \tau <0,
$$
and
we consider the equation
\begin{equation}
\label{forwardconparametri}
X_\tau = S(\tau-t)x+\int_0^\tau 1_{[t,T]}(\sigma) S(\tau-\sigma)
F(\sigma,X_\sigma)\; d\sigma
+\int_0^\tau 1_{[t,T]}(\sigma) S(\tau-\sigma)
G(\sigma,X_\sigma)\; dW_\sigma,
\end{equation}
for the unknown process $X_\tau$, $\tau\in [0,T]$.
Under the assumptions of
Hypothesis \ref{ipotesiuno}, equation
(\ref{forwardconparametri}) has
a  unique solution
$X\in \m S^2(H)$ which belongs to $\m S^p(H)$ 
for every $p\in [1,\infty)$.
  It clearly
 satisfies $X_\tau=x$ for $\tau\in [0,t)$, and
its restriction  to the time interval $[t,T]$
is the unique solution of (\ref{equazioneforwardmild}).
From now on we denote by $X(\tau ,t,x)$, $\tau\in [0,T]$,
the solution of (\ref{forwardconparametri}).

\begin{proposition}\label{xdifferenziabile}
  Assume Hypothesis \ref{ipotesiuno}. Then,
  for every $p\geq 1$,  the following hold.

  \begin{description}
    \item[(i)]   The map
  $(t,x)\mapsto X(\cdot,t,x)$ belongs to $\calg^{0,1}\Big([0,T]\times
  H, \m S^p(H)\Big)$.

    \item[(ii)]  For every $h\in H$, the directional derivative process $\nabla_xX(\tau,t,x)h$,
 $\tau\in [0,T]$, solves the equation:
$$
\left\{
\begin{array}{lll}\dis
\nabla_xX(\tau,t,x)h&=&\dis
 e^{(\tau-t)A}h+\int_t^\tau e^{(\tau-\sigma)A}
\nabla_xF(\sigma,X(\sigma,t,x))\nabla_xX(\sigma,t,x)h\; d\sigma
\\
&&+\dis
\int_t^\tau \nabla_x( e^{(\tau-\sigma)A}
G(\sigma,X(\sigma,t,x)))\nabla_xX(\sigma,t,x)h\; dW_\sigma,
\quad \tau\in [t,T],
\\\dis
\nabla_xX(\tau,t,x)h&=&h,\quad \tau\in [0,t).
\end{array}\right.
$$

\item[(iii)]  Finally
$\left\|\nabla_xX(\tau,t,x)h \right\|_{\infty, p}\leq c\,|h|$ for some constant $c$.
   \end{description}

\end{proposition}

\subsection{The backward equation}
\label{sec-backfor}
With this result in hands, we are know in position to study the regularity with respect to $(t,x)$ of the solution $(Y(\tau,t,x),Z(\tau,t,x))$ to the BSDE
\begin{equation}
\label{backwardforward}
 Y_\tau+\int_\tau^TZ_\sigma dW_\sigma= \phi(X(T,t,x))
  +\int_\tau^T\psi (\sigma ,X(\sigma,t,x), Y_\sigma,Z_\sigma)\, d\sigma.
\end{equation}
We assume that
$\psi: [0,T]\times H\times K \times L_2(\Xi,K)\to K$
verifies Hypothesis \ref{ipsupsi}.
 On the function $\phi :H\to K$ we make the following assumptions:
\begin{hypothesis}
\label{ipotesisuphi}
\begin{description}
\item[(i)] $\phi \in \calg^1 (H,K)$.
\item[(ii)]
  There exists $C>0$ and $m\geq 0$ such that,  
  $$
\forall x\in H,\quad \forall h \in H,\qquad   |\nabla_x\phi(x) h| \leq C\,|h| \left(1+|x|^m \right).
  $$
\end{description}
\end{hypothesis}

Under the assumptions of Hypotheses \ref{ipotesiuno} and
\ref{ipsupsi}, \ref{ipotesisuphi} by Propositions \ref{esisteunicox} and
 \ref{solbak} there exists a unique solution of
(\ref{backwardforward}) that we denote by  $(Y(\tau,t,x), Z(\tau,t,x))$, $\tau\in [0,T]$. This solution belongs to $\m K^p$ for each $p>1$.

Let us recall that the $X(\cdot,t,x)$ 
$\calf_{[t,T]}$-measurable, so that $Y(t,t,x)$ is measurable  with
respect to $\calf_{[t,T]}$ and $\calf_{t}$; it follows that $Y(t,t,x)$
is indeed deterministic (see also \cite{ElK}).

For later use we notice two useful identities:
 for $t\leq s \leq T$  the equality: $\P$-a.s.,
$$
X(\tau,s,X(s,t,x))=X(\tau,t,x), \qquad \tau\in [s,T],
$$
 is a  consequence of the uniqueness of the
solution of (\ref{equazioneforwardmild}). Since the solution
of the backward equation  is uniquely
determined on an interval $[s,T]$
by the values of the process $X$ on the same interval,
 for $t\leq s \leq T$ we have,
 $\P$-a.s.,
\begin{equation}
\label{markov}
Y(\tau,s,X(s,t,x))=Y(\tau,t,x),
\;\;  {\rm  for\; }\tau\in [ s, T],
\end{equation}
together with
$$
Z(\tau,s,X(s,t,x))=Z(\tau,t,x)
\;\; {\rm  for\; a.e.\; }\tau\in [ s, T].
$$

To investigate regularity properties of the dependence
on $t$ and $x$, we notice that the solution $(Y(\sigma,t,x),Z(\sigma,t,x))$ of the backward equation in (\ref{backwardforward}) can be written, with the notation of Propositions \ref{xdifferenziabile} and
\ref{solbak}, as 
$$
\left(Y(\cdot,t,x),Z(\cdot,t,x)\right)= \Phi\left(X(\cdot,t,x), \phi(X(T,t,x))\right).$$
Moreover, as a consequence of
Hypothesis \ref{ipotesisuphi}, it can be easily proved that $\xi\mapsto\phi(\xi)$ belongs to the space
$\calg^1(L^{(m+1)p}(\Omega ; H),L^p(\Omega ; K))$, for every
$p\in [1,\infty)$.
The following Proposition is then an immediate consequence of
Propositions  \ref{esisteunicox}, \ref{xdifferenziabile} and
 \ref{solbak},
and the chain rule for the class $\calg$, stated in
Lemma \ref{proprietadig}.

\begin{proposition}\label{ydifferenziabilenelsistema}
Assume Hypotheses \ref{ipotesiuno}, \ref{ipsupsi} and \ref{ipotesisuphi}.

Then the map
  $(t,x)\mapsto (Y(\cdot,t,x),Z(\cdot,t,x))$
  belongs to
$\calg^{0,1}([0,T]\times H\; ,\;\calk_{p})$ for all
$p\in (1,\infty)$.

  Denoting by $\nabla_xY$, $\nabla_xZ$
  the partial G\^ateaux derivatives with respect to $x$,
  the directional derivative process in the  direction $h\in H$,
$\left\{ \nabla_xY(\tau,t,x)h, \nabla_xZ(\tau,t,x)h \right\}_{\tau\in [0,T]}$ solves the
equation: $\P$-a.s., for $\tau\in[0,T]$, 
\begin{eqnarray*}
\lefteqn{\nabla_xY(\tau,t,x)h +\int_\tau^T\nabla_xZ(\sigma,t,x)h
\; dW_\sigma } \hspace*{20mm}\\
& & =  \nabla \phi (X(T,t,x))\nabla_x X(T,t,x)h \\
& & + \int_\tau^T
\nabla_x \psi (\sigma ,X(\sigma,t,x), Y(\sigma,t,x),Z(\sigma,t,x))\nabla_x
X(\sigma,t,x)h\; d\sigma
\\
&& +\int_\tau^T
\nabla_y
\psi (\sigma ,X(\sigma,t,x), Y(\sigma,t,x),Z(\sigma,t,x))
\nabla_xY(\sigma,t,x)h\; d\sigma
\\
&& +\int_\tau^T
\nabla_z
\psi (\sigma ,X(\sigma,t,x), Y(\sigma,t,x),Z(\sigma,t,x))
\nabla_xZ(\sigma,t,x)h\; d\sigma
.
\end{eqnarray*}
Finally the following estimate holds, for each $p>1$ :
$$
\E\left[ \sup_{\tau\in[0,T]}|\nabla_{x}Y(\tau,t,x)h|^{p} + \left(\int_{0}^{T}|\nabla_{x}Z(\sigma ,t,x)h|^{2}\,d \sigma
\right)^{p/2}\right]^{1/p} \leq C |h|(1+|x|^{(m+1)(n\vee 2)})
$$
for a constant $C$ depending on $p$, $\mu$, $T$ and $L$.
\end{proposition}

\begin{proof}
The first two statements comes from the chain rule for the class $\m G$.
 The final estimate follows from \eqref{stimanablaYZ} applied with
 $$X=X(\cdot,t,x),\ N=\nabla_{x}X(\cdot,t,x)h,\ \xi=\phi(X(T,t,x)),\
\zeta=\nabla\phi(X(T,t,x))\nabla_{x}X(T,t,x)h,
$$
taking into account that by
 Propositions  \ref{esisteunicox}  and \ref{xdifferenziabile} we have
  $$
  \|N\|_{\infty,p}\leq c |h|,\quad
 \|X\|_{\infty,p(m+1)}\leq c\,(1+ |x|),
 $$
and, by Hypothesis \ref{ipotesisuphi}, we
 also obtain
 $ \|\xi \|_{p}\leq c(1+ |x|)^{m+1}$,
 $\|\zeta \|_{p}\leq c|h| (1+|x|)^m$
 for a suitable constant $c$. 
 \end{proof}

\section{Application to nonlinear PDEs}
\label{sec-nlpdes}
\subsection{The nonlinear Kolmogorov equation}
\label{ssec-nlke}
We denote by $\calb_p(H)$ the set of measurable functions
$\phi:H\to\R$ with polynomial growth, i.e. such that $\sup_{x\in
H}|\phi(x)|(1+|x|^a)^{-1}<\infty$ for some  $a>0$.

Let $X(\tau,t,x)$, $\tau\in [t,T]$,
 denote the solution of the stochastic equation
$$
X_\tau=e^{(\tau-t)A}x + \int_t^\tau
e^{(\tau-\sigma)A}F(\sigma,X_\sigma)\, d\sigma
+\int_t^\tau
e^{(\tau-\sigma)A}G(\sigma,X_\sigma)\, dW_\sigma,
$$
where $A$, $F$, $G$, satisfy the assumptions in Hypothesis
\ref{ipotesiuno}.
The transition semigroup $P_{t,\tau}$ is defined
for arbitrary $\phi\in\calb_p(H)$ by the formula
$$
P_{t,\tau}[\phi](x)=\E\left[ \phi(X(\tau,t,x)) \right],\qquad x\in H.
$$
The estimate $
  \E\sup_{\tau\in[t,T]}|X_\tau |^p\leq C(1+|x|)^p$,
see (\ref{stimadeimomentidix}), shows that
$P_{t,\tau}$ is well defined as a linear operator from $
\calb_p(H)$ into itself; the semigroup property
$P_{t,s}P_{s,\tau}=P_{t,\tau}$, $t\leq s\leq \tau$,
is well known.
Let us denote by $\call_t$ the generator of $P_{t,\tau}$:
$$
\call_t[\phi](x)=\frac{1}{2}{\rm Trace }\left(
G(t,x)G(t,x)^*\nabla^2\phi(x)\right) + \< Ax+F(t,x),\nabla\phi(x)\>,
$$
where $\nabla\phi$ and $\nabla^2\phi$
are the first and the second G\^ateaux derivatives of
$\phi$ (identified with elements of $H$ and $L(H)$ respectively).
This definition is formal, since the domain of $\call_t$ is
not specified; however, if $\phi$ is sufficiently regular,
the function
$v(t,x)=P_{t,T}[\phi](x)$, is a classical solution of the backward
Kolmogorov equation:
$$
\partial_t v(t,x)
+\call_t [v(t,\cdot)](x) =0,\quad t\in [0,T],\,
x\in H,\qquad v(T,x)=\phi(x), \quad x\in H.
$$
We refer to \cite{DaZa1} and \cite{Za} for a detailed
exposition.
When $\phi$ is not regular,
the function $v$ defined by the formula 
$v(t,x)=P_{t,T}[\phi](x)$ can be considered as a
generalized solution of the
 backward
Kolmogorov equation.

Here we are interested in a generalization of this equation,
written formally as
\begin{equation}
\label{kolmogorovnonlineare}
  \left\{\begin{array}{l}\dis
\partial_t u(t,x) +\call_t [u(t,\cdot)](x) + 
\psi (t, x,u(t,x),G(t,x)^*\nabla_xu(t,x)) = 0,\quad t\in [0,T],\,
x\in H,\\
\dis u(T,x)=\phi(x).
\end{array}\right.
\end{equation}
We will refer to this equation as the nonlinear Kolmogorov
equation. $\psi: [0,T]\times H\times \R\times \Xi\to \R$ is a given
function verifying (\ref{ipsupsi}) and $\nabla_xu(t,x)$ is the
G\^ateaux derivative of
$u(t,x)$ with respect to $x$: it is identified with an element
of $H$, so that $G(t,x)^*\nabla_xu(t,x)\in \Xi$.

Now we define the notion of solution of the nonlinear Kolmogorov
equation. We consider  the variation of
constants formula for
 (\ref{kolmogorovnonlineare}):
\begin{equation}
\label{solmild}
  u(t,x) =\int_t^TP_{t,\tau}[
\psi (\tau, \cdot,u(\tau,\cdot),G(\tau,\cdot)^*
\nabla_xu(\tau,\cdot))
](x)\; d\tau+ P_{t,T}[\phi](x),\quad t\in [0,T],\,
x\in H,
\end{equation}
and we  notice that this formula  is meaningful,
provided $\psi(t,\cdot,\cdot,\cdot)$,
$u(t,\cdot)$ and $\nabla_xu(t,\cdot)$ have polynomial
growth.
We use this formula as a definition for the solution of
(\ref{kolmogorovnonlineare}):
\begin{definition}\label{defdisoluzionemild}
We say that a function
$u:[0,T]\times H\to\R$ is a mild solution of the nonlinear
Kolmogorov equation
(\ref{kolmogorovnonlineare}) if the following conditions hold:
\begin{description}
  \item[(i)]
$u\in\calg^{0,1}([0,T]\times H,\R)$;
  \item[(ii)]  there exists $C>0$ and $d\in \N$ such that
  $|\nabla_xu(t,x)h|\leq C|h|(1+|x|^{d})$ for all $t\in [0,T]$, $x\in H$,
$h\in H$;

  \item[(iii)] equality (\ref{solmild}) holds.
\end{description}
\end{definition}

\begin{remark}
An equivalent formulation
of (\ref{kolmogorovnonlineare})
or (\ref{solmild})
would be the following: we  consider
the G\^ateaux derivative $\nabla_xu(t,x)$ as an element of
$\Xi^*= L(\Xi,\R)=L_2(\Xi,\R)$, we  take a function
$\psi: [0,T]\times H\times \R\times L_2(\Xi,\R)\to \R$ and we
write the equation in the form
$$
\partial_t u(t,x)+\call_t [u(t,\cdot)](x) +
\psi (t, x,u(t,x),\nabla_xu(t,x)G(t,x)) =0.
$$
The two forms are clearly equivalent provided we identify
$\Xi^*= L_2(\Xi,\R)$ with $\Xi$ by the Riesz isometry. This will be done
in the sequel. In particular, although we keep the notation
in (\ref{kolmogorovnonlineare}),
we will sometimes consider $\psi$ as a real valued function
defined on $[0,T]\times H\times \R\times L_2(\Xi,\R)$,
satisfying Hypothesis \ref{ipsupsi} with
$K=\R$.
\end{remark}

We are now ready to state the main result of this section.
\begin{theorem}\label{main}
Assume that Hypothesis \ref{ipotesiuno}, \ref{ipsupsi}
(with $K=\R$) and
\ref{ipotesisuphi} hold.

The nonlinear Kolmogorov equation (\ref{kolmogorovnonlineare}) has a unique mild solution $u$ given by the formula
$$
u(t,x) =Y(t,t,x),\quad (t,x)\in[0,T]\times H
$$
where $(X,Y,Z)$ is the solution of the backward-forward
system
(\ref{backwardforward}). Moreover, we have, $\P$--a.s.
$$
Y(s,t,x) = u(s,X(s,t,x)),\qquad Z(s,t,x)= 
G(s, X(s,t,x))^*\nabla_xu(s,t, X(s,t,x)).
$$
\end{theorem}

\begin{proof}
Let us first recall a result of \cite[Lemma 6.3]{FT1}. 
Let $\{ e_i\}$ be a basis of  $\Xi$ and let us consider
the standard real Wiener process $W^i_\tau = \int_0^\tau \< e_i,dW_\sigma\>$, $\tau\geq 0$.

If $v\in \m G^{0,1}([0,T]\times H, \rset)$, for every $i$, the  quadratic variation 
  of $ u(s,X(s,t,x))$ and $W^i_s$ is given by
\begin{equation}
\label{quadrvardiu}
\left[v(\cdot,X(\cdot,t,x),W^i\right]_s = \int_t^s \nabla_x v(\tau, X(\tau,t,x))G(\tau, X(\tau,t,x))e_i
\; d\tau,\quad s\in [t,T].
\end{equation}

{\em (a) Existence.}
As we pointed out before, for $s\in[t,T]$, $Y(s,t,x)$ is $\m F_{[t,s]}$--measurable so that $Y(t,t,x)$ is deterministic. Moreover, as a byproduct of Proposition~\ref{ydifferenziabilenelsistema}, the function $u$ defined by the formula $u(t,x)=Y(t,t,x)$ has the
regularity properties stated in Definition
\ref{defdisoluzionemild}. It
remains to verify that
equality (\ref{solmild}) holds true for $u$.

To this purpose we first fix $t\in [0,T]$ and $x\in H$. Since $(Y(\cdot,t,x),Z(\cdot,t,x)$ solves the BSDE~\eqref{backwardforward}), we have, for $s\in[t,T]$,
$$
 Y(s,t,x)+\int_s^TZ(\tau,t,x)\, dW_\tau= \phi(X(T,t,x)) + 
  \int_s^T\psi \Big(\tau ,X(\tau,t,x),
   Y(\tau,t,x),Z(\tau,t,x)\Big)\, d\tau,
$$
and, taking expectation for $s=t$ we obtain, coming back to the definition of $u$ and $P_{t,T}$,
\begin{equation}
\label{presque}
 u(t,x)= P_{t,T}[\phi](x) + 
  \E\left[\int_t^T\psi \Big(\tau,X(\tau,t,x), Y(\tau,t,x),Z(\tau,t,x)
  \Big)\; d\tau \right].
\end{equation}
Moreover, we have, for each $i$,
$$
\left[ Y(\cdot,t,x),W^i\right]_s = \int_t^s \<Z_\tau, e_i \> \,d\tau, \quad s\in[t,T].
$$

Now let us observe that from the Markov property~\eqref{markov} we have, $\P$--a.s.,
$$
u(\t,X(\t,t,x))= Y(\t,t,x), \quad \t \in [t,T].
$$
It follows from \eqref{quadrvardiu} that, for each $i$,
$$
\left[ Y(\cdot,t,x),W^i\right]_s = \int_t^s \nabla_x u(\tau, X(\tau,t,x))G(\tau, X(\tau,t,x))e_i\, d\tau,\quad s\in [t,T].
$$
Therefore, for a.a. $\tau\in [t,T]$, we have $\P$-a.s.
$$
\nabla_xu(\tau, X(\tau,t,x))G(\tau, X(\tau,t,x))e_i=
\<Z(\tau,t,x),e_i\>,
$$
for every $i$. Identifying $\nabla_x u(t,x)$ with an element of
$\Xi$, we conclude that for a.a. $\tau\in [t,T]$,
$$
G(\tau, X(\tau,t,x))^*\nabla_xu(\tau,t, X(\tau,t,x))=
Z(\tau,t,x).
$$

Thus, $ \psi \left(\tau,X(\tau,t,x), Y(\tau,t,x),Z(\tau,t,x)
  \right) $ can be rewritten as
$$
\psi \left(\tau, X(\tau,t,x),u(\tau,X(\tau,t,x)),G(\tau,X(\tau,t,x))^*
\nabla_xu(\tau,X(\tau,t,x))\right)
$$
and \eqref{presque} leads to
$$
u(t,x)= P_{t,T}[\phi](x) + \int_t^TP_{t,\tau}[
\psi (\tau, \cdot,u(\tau,\cdot),G(\tau,\cdot)^*
\nabla_xu(\tau,\cdot))
](x)\, d\tau
$$
which is \eqref{solmild}.

\medskip

{\em (b) Uniqueness.}
Let $u$ be a mild solution. We look for a convenient
expression for the process
$u(s,X(s,t,x))$, $s\in [t,T]$.
By (\ref{solmild}) and the definition of $P_{t,\tau}$,
for every $s\in [t,T]$ and $x\in H$,
\begin{eqnarray*}
  u(s,x)  & = & \E \left[\phi(X(T,s,x)) \right] \\
 & & +\E\left[\int_s^T \psi\big(\tau, X(\tau,s,x),u(\tau,X(\tau,s,x)),
G(\tau,X(\tau,s,x))^*
\nabla_xu(\tau,X(\tau,s,x))\big)
d\tau\right].
\end{eqnarray*}

Since $X(\tau,s,x)$ is independent of $\calf_s$, we can replace the
expectation by the conditional expectation given $\calf_s$:
\begin{eqnarray*}
  u(s,x)  & = & \E^{\calf_s} \left[\phi(X(T,s,x)) \right] \\
 & & +\E^{\calf_s}\left[\int_s^T \psi\big(\tau, X(\tau,s,x),u(\tau,X(\tau,s,x)),
G(\tau,X(\tau,s,x))^*
\nabla_xu(\tau,X(\tau,s,x))\big)
d\tau\right].
\end{eqnarray*}

Taking into account the Markov property of $X$, $\P$--a.s.
$$
X(\tau,s,X(s,t,x))=X(\tau,t,x),
\qquad \tau\in [s,T],
$$
we have 
\begin{eqnarray*}
 \lefteqn{u(s,X(s,t,x)) =\E^{\calf_s} \left[ \phi(X(T,t,x)) \right]  } \\
   & & \qquad
+\E^{\calf_s}\left[\int_s^T
\psi\big(\tau, X(\tau,t,x),u(\tau,X(\tau,t,x)),
G(\tau,X(\tau,t,x))^*
\nabla_xu(\tau,X(\tau,t,x)) \big)
d\tau\right].
\end{eqnarray*}
If we set 
$$
\xi=\phi(X(T,t,x))+\int_t^T
\psi\big(\tau, X(\tau,t,x),u(\tau,X(\tau,t,x)),
G(\tau,X(\tau,t,x))^*
\nabla_xu(\tau,X(\tau,t,x))\big) d\tau
$$
the previous equality leads to
\begin{eqnarray*}
\lefteqn{u(s,X(s,t,x))} \\
&&  =\E^{\calf_s} \,[\xi] -\int_t^s
\psi \big(\tau, X(\tau,t,x),u(\tau,X(\tau,t,x)),
G(\tau,X(\tau,t,x))^*
\nabla_xu(\tau,X(\tau,t,x))\big)
 d\tau.
\end{eqnarray*}

Let us observe that $\E^{\calf_t}[\xi] = u(t,x)$.
Since $\xi\in L^2(\Omega;\R)$ is $\calf_{[t,T]}$--measurable,
by the representation theorem, there exists
$\widetilde{Z}\in L^2_\calp(\Omega \times [t,T]; L_2(\Xi,\R))$ such
that 
$$
\E^{\calf_s} [\xi] =u(t,x) + \int_t^s\widetilde{Z}_\tau\, dW_\tau,\quad s\in[t,T].
$$
We conclude that
 the process
$u(s,X(s,t,x))$, $s\in [t,T]$ is a (real) continuous
semimartingale with canonical decomposition
\begin{eqnarray}
\label{scomposizione}
u(s,X(s,t,x)) & = & u(t,x) + \int_t^s\widetilde{Z}_\tau\; dW_\tau\\
\nonumber
 &&-\int_t^s
\psi \Big(\tau, X(\tau,t,x),u(\tau,X(\tau,t,x)),
G(\tau,X(\tau,t,x))^*
\nabla_xu(\tau,X(\tau,t,x))\Big)
\, d\tau.
\end{eqnarray}

Using \eqref{quadrvardiu} and arguing as in the proof of existence, we deduce that
for a.a. $\tau\in [t,T]$, $\P$-a.s.
$$
G(\tau, X(\tau,t,x))^*\nabla_xu(\tau, X(\tau,t,x))=
\widetilde{Z}_\tau.
$$
Substituting into \eqref{scomposizione} we obtain
\begin{eqnarray*}
u(s,X(s,t,x)) & = &  u(t,x) +
\int_t^sG(\tau, X(\tau,t,x))^*\nabla_xu(\tau, X(\tau,t,x))\, dW_\tau \\
&& -\int_t^s
\psi \big(\tau, X(\tau,t,x),u(\tau,X(\tau,t,x)),
G(\tau,X(\tau,t,x))^*
\nabla_xu(\tau,X(\tau,t,x))\big)
\, d\tau,
\end{eqnarray*}
for $s\in [t,T]$. Since $u(T,X(T,t,x))=\phi(X(T,t,x))$, we deduce that
$$
\left\{ \big(u(s,X(s,t,x)), G(\tau,X(\tau,t,x))^*\nabla_x u(\tau,X(\tau,t,x))\big) \right\}_{s\in[t,T]}
$$
solves
the backward equation \eqref{backwardforward}. By uniqueness, we have $Y(s,t,x)=u(s,X(s,t,x))$, for each $s\in [t,T]$ and in particular, for $s=t$, $u(t,x)=Y(t,t,x)$.

\end{proof}

\subsection{Application to Optimal Control}
\label{ssec-control}
We wish to apply the above results to perform the synthesis of the
optimal control for a general nonlinear control system.

Fixed $t\in[0,T]$ and $x\in H$ an {\it admissible control
system} (a.c.s) is given
 by $(\Omega,\mathcal{E}, {\cal F}_{t}, \P, W_{t}, u)$ where
 \begin{itemize}
    \item  $(\Omega,\mathcal{E}, \P)$ is a probability  space,

    \item  $\{{\cal F}_{t}: t\geq 0\}$ is a filtration in it, satisfying
    the usual conditions,

    \item $\{W_{t}: t\geq  0\}$ is a cylindrical  $\P$-Wiener process
    with values in $\Xi$ and adapted to the filtration
    $\{{\cal F}_t\}$,

    \item  $u\in L_{\mathcal{P}}^{2}(\Omega\times [t,T]; U)$ satisfies the constraint:
 $u_t\in \mathcal{U}$ $\P$-a.s. for a.a.
 $t\in [t,T]$,
 where $U$ is a separable real Hilbert space and $\mathcal{U}$ is a fixed  subset of $U$.
 \end{itemize}
  To each a.c.s. we associate the mild solution
   $X^{u}\in C_{\mathcal{P}}([t,T];L^{2}(\Omega;H))$ of the state equation:
 \begin{equation}
 \label{equazionestato}
\left\{\begin{array}{l}\dis
dX^{u}_\tau =AX^{u}_\tau\,d\t + F(\tau,X^{u}_\tau)\,d\tau
+G(\tau,X^{u}_\tau)[r(\tau,X^{u}_{\tau},u_{\tau})\; d\tau
+ dW_\tau],\qquad \tau\in [t,T],\\
\dis X_{t}=x\in H.
\end{array}\right.
\end{equation}
and the cost:
\begin{eqnarray*}
\lefteqn{J(t,x,u) } \\ 
&= &  \E\left[ \int_{t}^{T}
\exp \left( \int_t^s \lambda(r,X_r^u,u_r)dr \right) l(s,X^{u}_s,u(s))\, ds + \exp \left( \int_t^T \lambda(r,X_r^u,u_r)dr\right)\phi (X^{u}_{T})\right]    
\end{eqnarray*}
where $\lambda, l$ and $\Phi$ are real functions. $\lambda$ may be called the discount function.

Our purpose is to minimize the functional $J$ over all a.c.s.
Notice that in this problem the discount $\lambda$ depends on the control $u$. We define the Hamiltonian function relative to the
above problem: for all $t\in [0,T]$, $x\in H$, $y\in \R$ and $z\in\Xi$
\begin{equation}
\psi(t,x,y,z)=\inf\left\{ l(t, x,u)+ <z,r(t,x,u)>_{\Xi} + \lambda(t,x,u)y : u\in
\mathcal{U}\right\}.
\label{definhamilton}
\end{equation}

We make the following assumption.
\begin{hypothesis}\label{ipotesicontrollo} The
following holds:
\begin{enumerate}
    \item $A, F$ and $G $ verify Hypothesis \ref{ipotesiuno}.

    \item $r:[0,T]\times H \times \mathcal{U} \rightarrow \Xi$ and there exists a constant $C>0$ such that
$$|r(t,x,u)|_{\Xi} \leq C$$

\item $l: [0,T] \times H \times \mathcal{U} \rightarrow \R$ satisfies the following condition
$$ 0 \leq l(t,0,u) \leq q_1(t) + h(u),$$
for some functions $ q_1 \in L^1(0,T)$ and $h:\m U \fl \rset_+$;
     \item $\lambda:[0,T]\times H \times \mathcal{U} \rightarrow \R$ is a non positive function;

     \item $\phi$ satisfies Hypothesis \ref{ipotesisuphi}.

      \item $\psi : [0,T] \times H \times \R \times \Xi \rightarrow \R$ is  measurable and, for $\sigma\in [0,T]$,
  $\psi (\sigma,\cdot,\cdot,\cdot)\in \mathcal{G}^{1,1,1} (H \times \R \times \Xi)$ with
  $$
  |\nabla_x\psi(\sigma,x, y,z)|+|\nabla_y\psi(\sigma,x, y,z)|\leq
  q(\sigma)+c\left(|x|^m+|y|^n+ |z|^2\right)
  $$
  where $q$ is in $\lp^1(0,T)$.

  \item For all $t \in [0,T]$, $x \in H$, $y \in \R$ and $z \in \Xi$ there exists a unique $\Gamma(t,x,y,z)\in\mathcal{U}$ that realizes the minimum in (\ref{definhamilton}). Namely
$$\psi(t,x,y,z)= l(t, x,\Gamma(t,x,y,z))+ <z,r(t,x,\Gamma(t,x,y,z))> + \lambda(t,x,\Gamma(t,x,y,z))y. $$
Moreover the map $\Gamma: [0,T] \times H \times \R \times \Xi \fl  U$ is measurable.
\end{enumerate}
\end{hypothesis}

\begin{remark}
The function
 $\psi: [t,T]\times H\times \R \times \Xi \rightarrow \R$ verifies
Hypothesis \ref{ipsupsi}. In particular, the Hypothesis \ref{ipotesicontrollo}-2. implies that $\psi$ is Lipschitz with respect to $z$.
The condition \ref{ipotesicontrollo}-3. implies the Hypothesis \ref{ipsupsi} (v). Finally, it follows from \ref{ipotesicontrollo}-4. that $\psi$ is nonincreasing so that we can take $\mu=0$ in Hypothesis \ref{ipsupsi} (ii).

Thus, if we identify  $\Xi$ with $ L_2(\Xi,\R)$,
the real function
 $\psi$ defined on $[t,T]\times H \times \R \times L_2(\Xi,\R) $ verifies
Hypothesis \ref{ipsupsi} with $K=\R$.
Therefore by Theorem \ref{main}  the Hamilton Jacobi Bellman
 equation 
 $$
  \left\{\begin{array}{l}\dis
\frac{\partial v(t,x)}{\partial t}+\call_t [v(t,\cdot)](x) +
\psi (t, x,v(t,x),G(t,x)^*\nabla_x v(t,x)) =0,\qquad t\in [0,T],\,
x\in H,\\
\dis v(T,x)=\phi(x).
\end{array}\right.
$$
admits a unique mild solution $v\in\m G^{0,1}([0,T]\times H)$.
\end{remark}

\begin{example}
Let us consider the following situation: $U=\rset$, $\mathcal{U}=\rset_+$ and $r(t,x,u)=0$, $l(t,x,u)=u^2/2$, $\l(t,x,u)=-u$. Then we have
$$
\psi(t,x,y,z)= -\frac{1}{2}\,y_+^2,\qquad \Gamma(t,x,y,z) = y_+.
$$
We see on this simple example that the Hamiltonian function $\psi$ is not Lipschitz with respect to $y$ but nonincreasing with a polynomial growth order.
\end{example}

We are in a position to prove the main result of this section:
 \begin{theorem}\label{maincontrollo}
 Assume Hypothesis \ref{ipotesicontrollo}.
 For all a.c.s. we have $J(t,x,u)\geq v(t,x)$ and the
 equality holds if and only if the following feedback law is verified
 by $u$ and $X^{u}$:
 \begin{equation}
  \label{leggefeedback} 
u(\sigma)= \Gamma(\sigma,X^{u}_{\sigma },v(\sigma,X^{u}_{\sigma}),
G(\sigma,X^{u}_{\sigma})\nabla_{x}v(\sigma,X^{u}_{\sigma})),\quad
\P- {\rm a.s. \; for\;a.a.\; } \sigma\in [t,T].
   \end{equation}
    Finally there exists at least an a.c.s. for which
    (\ref{leggefeedback}) holds. In such a system the closed loop
    equation:
     \begin{equation}
     \label{equazioneclosedloop}
\left\{\begin{array}{l}
\dis{d \overline{X}_\tau
= A\overline{X}_{\tau}\,d\t +F(\tau,\overline{X}_\tau)\, d\tau+G(\tau,\overline{X}_\tau)\, dW_\tau  }\\ 
\qquad\dis{+G(\tau,\overline{X}_\tau)r(\tau,\overline{X}_{\tau},\Gamma(\tau,
\overline{X}_{\tau},v(\t,\overline{X}_{\t}),
G(\t,\overline{X}_{\t})\nabla_{x}v(\t,\overline{X}_{\t})))\, d\tau }\\
\\

 \dis \overline{X}_{t}=\dis x\in H.
\end{array}\right.
\end{equation}
admits a solution and if
$\overline{u}(\sigma)= \Gamma(\sigma,\overline{X}_{\sigma}
 ,v(\sigma,\overline{X}_{\sigma}),
G(\sigma,\overline{X}_{\sigma})\nabla_{x}v(\sigma,\overline{X}_{\sigma}))$  then the couple
 $(\overline{u},\overline{X})$ is optimal for the control problem.
 \end{theorem}

\begin{proof}
For all a.c.s.,
setting $u(s)=0$ for
$s<t$,  the Girsanov theorem ensures that there
exists a probability measure $\widetilde{\P}$ on  $\Omega$ such that
$$\widetilde{W}_{t}:=
W_{t}+\int_{0}^{t}r(\sigma,X^{u}_\sigma ,u_{\sigma} )
\; d\sigma $$
is a $\widetilde{\P}$-wiener process
(notice that the function $r$ is bounded). Relatively to $\widetilde{W}$ equation (\ref{equazionestato})
can be rewritten:
 $$
dX^{u}_\tau =AX^{u}_{\tau}\, d\t + F(\tau,X^{u}_\tau)\, d\tau
+G(\tau,X^u_{\tau})\, d\widetilde{W}_\tau,\quad \tau\in [t,T], \qquad X_{t}^u=x\in H.
$$
The process
$X^{u}$ turns out to be adapted to the filtration
$\widetilde{\calf}_{t}$ generated by $\widetilde{W}$ and
completed in the usual way.
 In the space $(\Omega, {\mathcal E},
\{\widetilde{\mathcal{F}}_{t}\},\widetilde{\P})$,
we can solve the SDE
$$
\widetilde{X}(\tau,t,x) = e^{(\tau-t)A}x+\int_t^\tau e^{(\tau-\sigma)A}
F(\sigma,\widetilde{X}(\sigma,t,x))\; d\sigma
+\int_t^\tau e^{(\tau-\sigma)A}
G(\sigma,\widetilde{X}(\sigma,t,x))\; d\widetilde{W}_\sigma
$$
and then the BSDE
$$
\widetilde{Y}(\tau,t,x)= \phi(\widetilde{X}(T,t,x))
  +\int_\tau^T\psi (\sigma
,\widetilde{X}(\sigma,t,x),\ty(\s,t,x),\widetilde{Z}(\sigma,t,x)) d\sigma-
 \int_\tau^T\widetilde{Z}(\sigma,t,x)
d\widetilde{W}_\sigma  .
$$
Actually, this construction can be done for arbitrary $t\in[0,T]$ and $x\in H$. Let us set 
$$
\forall \tau\in[t,T],\qquad D_\tau = \exp\left(\int_t^\t
\lambda\left(r,\widetilde{X}(r,t,x),u_r\right)dr\right).
$$
We have
\begin{eqnarray*}
d\, D_\tau \ty(\t,t,x)  & = & D_\tau\left\{\lambda(\t,\widetilde{X}(\t,t,x),u_\t)\ty(\t,t,x) -
\psi(\t,\widetilde{X}(\t,t,x),\widetilde{Y}(\t,t,x),\widetilde{Z}(\t,t,x)) \right\} d\t \\
 & &  + D_\tau \widetilde{ Z}(\t,t,x) d\widetilde{W}_\t 
 \end{eqnarray*}
and thus
\begin{eqnarray*}
D_T\phi(X_T) & = & \widetilde{Y}(t,t,x) + \int_t^T D_\tau\widetilde{Z}(\t,t,x) d\widetilde{W}_\t \\
  & &+  \int_t^T D_\tau\{\lambda(\t,\widetilde{X}(\t,t,x),u_\t)\widetilde{Y}(\t,t,x)
 - \psi(\t,\widetilde{X}(\t,t,x),\widetilde{Y}(\t,t,x),\widetilde{Z}(\t,t,x))\} d
\t 
\end{eqnarray*}

We notice that $\widetilde{X}(\sigma,t,x)=X^{u}_\sigma$ and coming back to the original Wiener process we get:
 \begin{eqnarray*}
  D_T\phi(X_T^u) & = & \ty(t,t,x) +  \int_{t}^T D_\tau \left\{ \tz(\t,t,x) r(\t,X_{\t}^u,u_{\t}) d\t +
   \tz(\t,t,x) dW_\t \right\}\\
  & & + \int_{t}^T D_\tau \left\{\lambda(\t,X_\t^u,u_\t)\ty(\t,t,x)
   - \psi(\t,X_\t^u,\ty(\t,t,x),\tz(\t,t,x))\right\}.
   \end{eqnarray*}
Now, from Proposition~\ref{main}, we have  $
    \widetilde{Y}(\tau,t,x)=v(t,X^u_\tau)$ and 
    $$
      \widetilde{Z}(\tau,t,x)=G(\tau,
\widetilde{X}(\tau,t,x))^{*}\nabla_x v(\tau,
\widetilde{X}(\tau,t,x))
=G(\tau,
X^u_\tau)^{*}\nabla_x v(\tau,
X^u_\tau).
   $$
Taking expectation  with respect to
 the original probability $\P$ in
 the previous relation we obtain:
\begin{eqnarray*}
\lefteqn{ \E\, \left[  D_T
\phi(X^{u}_T)\right] -v(t,x) }\\
& = &\E\left[ \int_{t}^T D_\sigma \left\{\lambda(\s,X_\s^u,u_\s)v(\s,X_\s^u)
  -\psi(\sigma ,X^{u}_\sigma,v(\s,X_\s^u),
 G(\sigma ,X^{u}_\sigma)^*\nabla_{x}v(\sigma ,X^{u}_\sigma))
\,  d\sigma \right\} \right] \\
 & & +  \E\left[ \int_{t}^T D_\sigma G(\sigma ,X^{u}_\sigma)^*\nabla_{x}v(\sigma ,X^{u}_\sigma))r(\s,X_{\s}^u,u_{\s})
  \, d\sigma\right].
  \end{eqnarray*}
  Adding to both hand side the term
 $$
 \E\left[ \int_{t}^{T}  D_\sigma l(\sigma, X^{u}_{\sigma},u_{\sigma })\,d \sigma \right]
 $$
we get to the following expression
$$
 J(t,x,u) - v(t,x)  =    \E\left[ \int_{t}^T D_\sigma\, H\left(\sigma,X^u_\sigma,v(\sigma ,X^{u}_\sigma), G(\sigma ,X^{u}_\sigma)^*\nabla_{x}v(\sigma ,X^{u}_\sigma),u_\sigma \right)
\,  d\sigma \right], 
$$
where we have set
$$
H(\sigma,x,y,z,u)=-\psi(\sigma ,x,y,z)  +   \lambda(\s,x,u)y  +   z
r(\s,x,u)+ l(\sigma, x,u)
$$
 The above equality is known as the {\it fundamental relation} and
 immediately implies, by definition of $\psi$  that $ v(t,x)\leq J(t,x,u)$
  and that the equality holds if and only if (\ref{leggefeedback}) holds.

 Finally the existence of a weak solution to equation
(\ref{equazioneclosedloop})
 is again a consequence of the Girsanov theorem. Namely let $X\in
 C_{\mathcal{P}}([t,T];L^{2}(\Omega;H))$ be the mild solution of
 $$\left\{
 \begin{array}{l}
 {\displaystyle dX_{\tau}=AX_{\tau}\,d\t + F(\tau,X_{\tau})\,d\tau
 +G(\tau,X_{\tau})\,dW_{\tau} }, \\
    X_{t}=x,
 \end{array}
 \right.
 $$ 
 and let $\widehat{P}$ be the probability on $\Omega$ under
 which
 $$\widehat{W}_{t}:= -\int_{0}^{t} r(\sigma
  ,X_{\sigma },\Gamma(\sigma,X_{\sigma },r(\sigma
  ,X_{\sigma },G(\sigma,X_{\sigma})^*
  \nabla_{x}v(\sigma,X^{u}_{\sigma}))\;
  d\sigma+W_{t}$$
is a Wiener process ($r$ is bounded). Then $X$ is the mild solution of equation
(\ref{equazioneclosedloop}) relatively to the probability $\widehat{P}$
and the Wiener process $\widehat{W}$. 
\end{proof}


\begin{thebibliography}{200}

\bibitem{Be} Bensoussan, A.; \textit{Stochastic control by functional analysis methods.} 
Studies in Mathematics and its Applications, 11. 
North-Holland Publishing Co., Amsterdam-New York, 1982.


\bibitem{BeLi} Bensoussan, A.; Lions, J.L.;
\textit{Application of variational inequalities in stochastic control.}
Studies in Mathematics and its Applications, 12. 
North-Holland Publishing Co., Amsterdam-New York, 1982. 

\bibitem{BDHPS} Briand, Ph.; Delyon, B.; Hu, Y.; Pardoux, E.; Stoica, L.
\textit{$L\sp p$ solutions of backward stochastic differential equations.} Stochastic
Process. Appl. 108 (2003), no. 1, 109-129.



\bibitem{C} Cerrai, S.;
\textit{Second order PDE's in finite and infinite dimensions. A probabilistic approach.} Probab. Theory Relat. Fields,115,1999,383-399.

\bibitem{DaZa1} G. Da Prato, J. Zabczyk,
{\bf Stochastic equations in infinite dimensions},
Encyclopedia of Mathematics and its Applications, 44,
Cambridge University Press, 1992.


\bibitem{ElK} N. El Karoui,
{\em Backward stochastic differential equations a general introduction},
 in: {\bf Backward Stochastic Differential Equations}, ed.  N. El
 Karoui, L. Mazliak 7-26. Pitman Research Notes in Mathematics Series
 364, Longman, 1997.

\bibitem{FT1} Fuhrman, M.; Tessitore, G.; \textit{Non linear Kolmogorov
equations in infinite dimensional spaces: the backward stochastic differential equations
approach and applications to optimal control. }Ann. Probab. 30 (2002), no. 3,
1397--1465.

\bibitem{FT2} Fuhrman, M.; Tessitore, G.; \textit{Infinite horizon backward stochastic
differential equations and elliptic equations in Hilbert spaces.} Ann. Probab. 32
(2004), no. 1B, 607--660.


\bibitem{FlSo} W. H. Fleming, H. M.  Soner,
{\bf Controlled Markov processes and viscosity solutions},
Springer-Verlag, 1993.



\bibitem{ElPQ}  N. El Karoui; S. Peng; M. C. Quenez, \textit{Backward stochastic differential equations in finance. Math. Finance} 7 (1997), no. 1, 1--71. 



\bibitem{GoRo} F. Gozzi, E. Rouy, \textit{Regular solutions of second-order stationary Hamilton-Jacobi equations}. (English. English summary) 
J. Differential Equations 130 (1996), no. 1, 201--234.



\bibitem{MaYo} J. Ma, J. Yong,
 {\bf Forward-backward stochastic differential
 equations and their applications},
Lecture Notes in Mathematics 1702,
Springer, 1999.


\bibitem{P2}Pardoux, E.;
\textit{Backward stochastic differential equations and viscosity
solutions of systems of semilinear parabolic and elliptic PDEs of
second order}. Stochastic analysis and related topics, VI (Geilo,
1996), 79--127, Progr. Probab., 42, Birkhäuser Boston, Boston, MA,
1998.

\bibitem{P3} Pardoux, E.; \textit{BSDEs, weak convergence and homogenization
of semilinear PDEs.} Nonlinear analysis, differential equations
and control (Montreal, QC, 1998), 503--549, NATO Sci. Ser. C Math.
Phys. Sci., 528, Kluwer Acad. Publ., Dordrecht, 1999.






\bibitem{PaPe} E. Pardoux, S. Peng,
{\em Adapted solution of a backward stochastic
differential equation}, Systems and Control Lett.
{\bf 14}, 1990, 55-61.


\bibitem{PaPe2} E. Pardoux, S. Peng,
{\em Backward stochastic differential equations and
quasilinear parabolic partial differential equations},
in:  {\bf Stochastic partial differential equations
and their applications}, eds. B.L. Rozowskii, R.B.
Sowers, 200-217,
Lecture Notes in Control Inf. Sci. 176, Springer, 1992.

\bibitem{Pe1} Peng, S.
\textit{Probabilistic interpretation for Systems of Quasilinear Parabolic Partial Differential Equations}
Stochastics, 1991, 37, 61-74. 

\bibitem{Pe2} Peng, S.
\textit{Stochastic Hamilton-Jacobi-Bellman equations.} SIAM J. Control Optim., 30, 284-304,
  1992.

\bibitem{Pe3} Peng, S.
\textit{A generalized dynamic programming principle and Hamilton-Jacobi-Bellman Equation}
Stochastics, 38, 119-134, 1992.

\bibitem{Pe4} Peng, S.
\textit{A nonlinear  Feynman-Kac formula and applications} in Proceedings of Symposium of System Sciences and Control Theory, ed Chen and Yong. Singapore: Word Scientific, 173-184.



\bibitem{Za} J. Zabczyk,  {\em Parabolic equations
on Hilbert spaces}, in: {\bf Stochastic PDE's and Kolmogorov Equations
in Infinite Dimensions}, ed. G. Da Prato, 117-213.
Lecture Notes in Mathematics 1715,
Springer, 1999.

\end{thebibliography}
\end{document}